\documentclass[reqno,12pt,a4paper]{amsart}

\usepackage{mathrsfs}
\usepackage{amssymb}
\usepackage{amsfonts}
\usepackage{latexsym}
\usepackage{amsthm}

\usepackage{graphicx}
\def\lb{\label}

\newcommand{\er}[1]{\textrm{(\ref{#1})}}

\begin{document}


\renewcommand{\theequation}{\arabic{section}.\arabic{equation}}
\theoremstyle{plain}
\newtheorem{theorem}{\bf Theorem}[section]
\newtheorem{lemma}[theorem]{\bf Lemma}
\newtheorem{corollary}[theorem]{\bf Corollary}
\newtheorem{proposition}[theorem]{\bf Proposition}
\newtheorem{definition}[theorem]{\bf Definition}

\def\a{\alpha}  \def\cA{{\mathcal A}}     \def\bA{{\bf A}}  \def\mA{{\mathscr A}}
\def\b{\beta}   \def\cB{{\mathcal B}}     \def\bB{{\bf B}}  \def\mB{{\mathscr B}}
\def\g{\gamma}  \def\cC{{\mathcal C}}     \def\bC{{\bf C}}  \def\mC{{\mathscr C}}
\def\G{\Gamma}  \def\cD{{\mathcal D}}     \def\bD{{\bf D}}  \def\mD{{\mathscr D}}
\def\d{\delta}  \def\cE{{\mathcal E}}     \def\bE{{\bf E}}  \def\mE{{\mathscr E}}
\def\D{\Delta}  \def\cF{{\mathcal F}}     \def\bF{{\bf F}}  \def\mF{{\mathscr F}}
\def\c{\chi}    \def\cG{{\mathcal G}}     \def\bG{{\bf G}}  \def\mG{{\mathscr G}}
\def\z{\zeta}   \def\cH{{\mathcal H}}     \def\bH{{\bf H}}  \def\mH{{\mathscr H}}
\def\e{\eta}    \def\cI{{\mathcal I}}     \def\bI{{\bf I}}  \def\mI{{\mathscr I}}
\def\p{\psi}    \def\cJ{{\mathcal J}}     \def\bJ{{\bf J}}  \def\mJ{{\mathscr J}}
\def\vT{\Theta} \def\cK{{\mathcal K}}     \def\bK{{\bf K}}  \def\mK{{\mathscr K}}
\def\k{\kappa}  \def\cL{{\mathcal L}}     \def\bL{{\bf L}}  \def\mL{{\mathscr L}}
\def\l{\lambda} \def\cM{{\mathcal M}}     \def\bM{{\bf M}}  \def\mM{{\mathscr M}}
\def\L{\Lambda} \def\cN{{\mathcal N}}     \def\bN{{\bf N}}  \def\mN{{\mathscr N}}
\def\m{\mu}     \def\cO{{\mathcal O}}     \def\bO{{\bf O}}  \def\mO{{\mathscr O}}
\def\n{\nu}     \def\cP{{\mathcal P}}     \def\bP{{\bf P}}  \def\mP{{\mathscr P}}
\def\r{\rho}    \def\cQ{{\mathcal Q}}     \def\bQ{{\bf Q}}  \def\mQ{{\mathscr Q}}
\def\s{\sigma}  \def\cR{{\mathcal R}}     \def\bR{{\bf R}}  \def\mR{{\mathscr R}}
\def\vs{\varsigma}  \def\cS{{\mathcal S}}     \def\bS{{\bf S}}  \def\mS{{\mathscr S}}
\def\t{\tau}    \def\cT{{\mathcal T}}     \def\bT{{\bf T}}  \def\mT{{\mathscr T}}
\def\f{\phi}    \def\cU{{\mathcal U}}     \def\bU{{\bf U}}  \def\mU{{\mathscr U}}
\def\F{\Phi}    \def\cV{{\mathcal V}}     \def\bV{{\bf V}}  \def\mV{{\mathscr V}}
\def\P{\Psi}    \def\cW{{\mathcal W}}     \def\bW{{\bf W}}  \def\mW{{\mathscr W}}
\def\o{\omega}  \def\cX{{\mathcal X}}     \def\bX{{\bf X}}  \def\mX{{\mathscr X}}
\def\x{\xi}     \def\cY{{\mathcal Y}}     \def\bY{{\bf Y}}  \def\mY{{\mathscr Y}}
\def\X{\Xi}     \def\cZ{{\mathcal Z}}     \def\bZ{{\bf Z}}  \def\mZ{{\mathscr Z}}
\def\O{\Omega}


\def\mb{{\mathscr b}}
\def\mh{{\mathscr h}}
\def\me{{\mathscr e}}
\def\mk{{\mathscr k}}
\def\mz{{\mathscr z}}
\def\mx{{\mathscr x}}

\def\be{{\bf e}} \def\bc{{\bf c}} \def\bt{{\bf t}}
\def\bx{{\bf x}} \def\by{{\bf y}}
\def\bv{{\bf v}} 
\def\Om{\Omega} \def\bp{{\bf p}}

\newcommand{\gA}{\mathfrak{A}}          \newcommand{\ga}{\mathfrak{a}}
\newcommand{\gB}{\mathfrak{B}}          \newcommand{\gb}{\mathfrak{b}}
\newcommand{\gC}{\mathfrak{C}}          \newcommand{\gc}{\mathfrak{c}}
\newcommand{\gD}{\mathfrak{D}}          \newcommand{\gd}{\mathfrak{d}}
\newcommand{\gE}{\mathfrak{E}}
\newcommand{\gF}{\mathfrak{F}}           \newcommand{\gf}{\mathfrak{f}}
\newcommand{\gG}{\mathfrak{G}}           \newcommand{\gog}{\mathfrak{g}}
\newcommand{\gH}{\mathfrak{H}}           \newcommand{\gh}{\mathfrak{h}}
\newcommand{\gI}{\mathfrak{I}}           \newcommand{\gi}{\mathfrak{i}}
\newcommand{\gJ}{\mathfrak{J}}           \newcommand{\gj}{\mathfrak{j}}
\newcommand{\gK}{\mathfrak{K}}            \newcommand{\gk}{\mathfrak{k}}
\newcommand{\gL}{\mathfrak{L}}            \newcommand{\gl}{\mathfrak{l}}
\newcommand{\gM}{\mathfrak{M}}            \newcommand{\gm}{\mathfrak{m}}
\newcommand{\gN}{\mathfrak{N}}            \newcommand{\gn}{\mathfrak{n}}
\newcommand{\gO}{\mathfrak{O}}
\newcommand{\gP}{\mathfrak{P}}             \newcommand{\gp}{\mathfrak{p}}
\newcommand{\gQ}{\mathfrak{Q}}             \newcommand{\gq}{\mathfrak{q}}
\newcommand{\gR}{\mathfrak{R}}             \newcommand{\gr}{\mathfrak{r}}
\newcommand{\gS}{\mathfrak{S}}              \newcommand{\gs}{\mathfrak{s}}
\newcommand{\gT}{\mathfrak{T}}             \newcommand{\gt}{\mathfrak{t}}
\newcommand{\gU}{\mathfrak{U}}             \newcommand{\gu}{\mathfrak{u}}
\newcommand{\gV}{\mathfrak{V}}             \newcommand{\gv}{\mathfrak{v}}
\newcommand{\gW}{\mathfrak{W}}             \newcommand{\gw}{\mathfrak{w}}
\newcommand{\gX}{\mathfrak{X}}               \newcommand{\gx}{\mathfrak{x}}
\newcommand{\gY}{\mathfrak{Y}}              \newcommand{\gy}{\mathfrak{y}}
\newcommand{\gZ}{\mathfrak{Z}}             \newcommand{\gz}{\mathfrak{z}}

\def\ve{\varepsilon} \def\vt{\vartheta} \def\vp{\varphi}  \def\vk{\varkappa}
\def\vr{\varrho}

\def\A{{\mathbb A}} \def\B{{\mathbb B}} \def\C{{\mathbb C}}
\def\dD{{\mathbb D}} \def\E{{\mathbb E}} \def\dF{{\mathbb F}} \def\dG{{\mathbb G}}
\def\H{{\mathbb H}}\def\I{{\mathbb I}} \def\J{{\mathbb J}} \def\K{{\mathbb K}}
\def\dL{{\mathbb L}}\def\M{{\mathbb M}} \def\N{{\mathbb N}} \def\dO{{\mathbb O}}
\def\dP{{\mathbb P}} \def\dQ{{\mathbb Q}} \def\R{{\mathbb R}}\def\S{{\mathbb S}} \def\T{{\mathbb T}}
\def\U{{\mathbb U}} \def\V{{\mathbb V}}\def\W{{\mathbb W}} \def\X{{\mathbb X}}
\def\Y{{\mathbb Y}} \def\Z{{\mathbb Z}}

\def\dk{{\Bbbk}}


\def\la{\leftarrow}              \def\ra{\rightarrow}     \def\Ra{\Rightarrow}
\def\ua{\uparrow}                \def\da{\downarrow}
\def\lra{\leftrightarrow}        \def\Lra{\Leftrightarrow}
\newcommand{\abs}[1]{\lvert#1\rvert}
\newcommand{\br}[1]{\left(#1\right)}

\def\lan{\langle} \def\ran{\rangle}


\def\lt{\biggl}                  \def\rt{\biggr}
\def\ol{\overline}               \def\wt{\widetilde}
\def\no{\noindent}


\let\ge\geqslant                 \let\le\leqslant
\def\lan{\langle}                \def\ran{\rangle}
\def\/{\over}                    \def\iy{\infty}
\def\sm{\setminus}               \def\es{\emptyset}
\def\ss{\subset}                 \def\ts{\times}
\def\pa{\partial}                \def\os{\oplus}
\def\om{\ominus}                 \def\ev{\equiv}
\def\iint{\int\!\!\!\int}        \def\iintt{\mathop{\int\!\!\int\!\!\dots\!\!\int}\limits}
\def\wh{\widehat}
\def\bs{\backslash}
\def\na{\nabla}
\def\ti{\tilde}
\def\hb{\hbar}
\def\cd{\centerdot}
\def\we{\wedge}

\def\sh{\mathop{\mathrm{sh}}\nolimits} \def\ch{\mathop{\mathrm{ch}}\nolimits}
\def\all{\mathop{\mathrm{all}}\nolimits}
\def\Area{\mathop{\mathrm{Area}}\nolimits}
\def\arg{\mathop{\mathrm{arg}}\nolimits}
\def\const{\mathop{\mathrm{const}}\nolimits}
\def\det{\mathop{\mathrm{det}}\nolimits}
\def\diag{\mathop{\mathrm{diag}}\nolimits}
\def\diam{\mathop{\mathrm{diam}}\nolimits}
\def\dim{\mathop{\mathrm{dim}}\nolimits}
\def\dist{\mathop{\mathrm{dist}}\nolimits}
\def\Im{\mathop{\mathrm{Im}}\nolimits}
\def\Iso{\mathop{\mathrm{Iso}}\nolimits}
\def\Ker{\mathop{\mathrm{Ker}}\nolimits}
\def\Lip{\mathop{\mathrm{Lip}}\nolimits}
\def\rank{\mathop{\mathrm{rank}}\limits}
\def\Ran{\mathop{\mathrm{Ran}}\nolimits}
\def\Re{\mathop{\mathrm{Re}}\nolimits}
\def\Res{\mathop{\mathrm{Res}}\nolimits}
\def\res{\mathop{\mathrm{res}}\limits}
\def\sign{\mathop{\mathrm{sign}}\nolimits}
\def\span{\mathop{\mathrm{span}}\nolimits}
\def\supp{\mathop{\mathrm{supp}}\nolimits}
\def\Tr{\mathop{\mathrm{Tr}}\nolimits}
\def\BBox{\hspace{1mm}\vrule height6pt width5.5pt depth0pt \hspace{6pt}}
\def\where{\mathop{\mathrm{where}}\nolimits}
\def\as{\mathop{\mathrm{as}}\nolimits}
\def\1{1\!\!1}



\newcommand\nh[2]{\widehat{#1}\vphantom{#1}^{(#2)}}
\def\di{\diamond}

\def\Oplus{\bigoplus\nolimits}



\def\qqq{\qquad}
\def\qq{\quad}
\let\ge\geqslant
\let\le\leqslant
\let\geq\geqslant
\let\leq\leqslant
\newcommand{\ca}{\begin{cases}}
\newcommand{\ac}{\end{cases}}
\newcommand{\ma}{\begin{pmatrix}}
\newcommand{\am}{\end{pmatrix}}
\renewcommand{\[}{\begin{equation}}
\renewcommand{\]}{\end{equation}}
\def\eq{\begin{equation}}
\def\qe{\end{equation}}
\def\[{\begin{equation}}
\def\bu{\bullet}
\def\ced{\centerdot}
\def\tes{\textstyle}
\def\ci{\circ}


\title[{Inverse problems for ZS-operators and their isomorphisms
 }]
{Inverse problems for ZS-operators and their isomorphisms}

\date{\today}

\author[Evgeny Korotyaev]{Evgeny Korotyaev}
\author[Zongfeng Zhang]{Zongfeng Zhang}
\address{Academy for Advance interdisciplinary Studies, Northeast Normal University,
Changchun, China, \ korotyaev@gmail.com}


\begin{abstract}
\no Consider two inverse problems for ZS-operators  problems on the
 unit interval. It means that there  are two corresponding mappings
$F, f$ from a Hilbert space of potentials   $H$ into their spectral
data. They are called isomorphic if $F$ is a composition of $f$ and
some isomorphism  $U$ of $H$ onto itself.
We consider isomorphic inverse problems for ZS-operators  on the unit interval under basic boundary conditions and on the circle.
The proof is based on the non-linear analysis and properties of the 4-spectra mapping constructed in our paper.

\end{abstract}

\maketitle

\begin{quotation}
\begin{center}
{\bf Table of Contents}
\end{center}

\vskip 6pt

{\footnotesize

1. Introduction and main results \hfill \pageref{Sec1}\ \ \ \ \

2. Preliminary results.
 \hfill
\pageref{Sec2}\ \ \ \ \

3. Gauge transformations and their properties  \hfill \pageref{Sec3}\ \ \ \ \

4. Isomorphic inverse problems on the unit interval \hfill \pageref{Sec4}\ \
\ \ \

5. Isomorphic inverse problems on the circle and replacing mappings \hfill \pageref{Sec5}\ \ \ \ \



 }
\end{quotation}



\section {Introduction and main results \lb{Sec1}}
\setcounter{equation}{0}

\subsection{Introduction}

We consider inverse problems for Zakharov-Shabat systems (or shortly ZS-systems) on the unit interval and show that they are isomorphic.
There are a lot of results about inverse problems.
Different results and  approaches to inverse spectral problems can be found in the monographs \cite{AT04}, \cite{M86}, \cite{L87}, \cite{LS91} \cite{PT87},  and
references therein.
In general, the study of
inverse spectral problems consists of the following parts:

\noindent (i) Uniqueness: prove that the spectral data
(eigenvalues plus some additional parameters) determine the potential uniquely);\\
(ii) Reconstruction: reconstruct the potential from spectral data;\\
(iii) Characterization: describe all spectral data corresponding to
fixed classes of potentials. \\
(iv) Stability estimates: obtain a priori two sided estimates of the
potential and spectral data.


We will discuss their additional {\it isomorphic} properties.

\no {\bf Definition.} {\it Let $f$ and $f_o$ be mappings from a
Hilbert space ${\mathcal K}$ to a set $X$. They are called isomorphic if
$f_o=f\circ U$ for some isomorphism (in general, non-linear) $U$ of
$\mathcal K$ onto itself.}

Note that if some of two inverse problems is a bijection, then $U$ is a unique
canonical automorphism of $\mathcal K$.  We shortly describe properties of isomorphic inverse problems. Assume that we have two isomorphic inverse problems, then we have

\no 1) If the first one has some property from (i)-(iv), then the second
also has it. For example, the first has uniqueness iff the second
has uniqueness.

\no 2) Eigenvalues of the first problem  have some asymptotics for
each potential iff
 eigenvalues of the second  problem have similar asymptotics.

\no  3) The first problem  has some trace formula iff
 the second  problem  has a similar trace formula.

Recall that isomorphic inverse problems for
Sturm-Liouville  problems on the unit interval and the circle
were described by Korotyaev \cite{K25} and we will use these results.

There are a lot of results about the inverse
problems for ZS-systems (or Dirac systems), see \cite{LS91},  \cite{A14}, \cite{AHM05}, \cite{DK00},  \cite{GD75} \cite{S70} on the unit interval
 under boundary condiotions and, see \cite{LS91}, \cite{BG93}, \cite{GK14}, \cite{K01}, \cite{K05}, \cite{K96} on the circle, and references therein.
We consider the ZS-systems on the interval $[0,1]$ under Dirichlet and Neumann
  boundary conditions
\begin{equation}
\label{DN} Jf'+Vf=\l f,\qqq
\begin{aligned}
&     \qqq f_1(0)= \textstyle f_1(1)=0, \qq \{\m_n, n\in \Z\} \ Dirichlet
\\
&   \qqq  f_2(0)=\textstyle f_2(1)=0,  \qq \{\n_n, n\in \Z\} \ Neumann
\end{aligned}\ \ \ ,
\end{equation}
where $f=(f_1,f_2)^\top$ is the vector function
and under the so-called mixed boundary conditions:
\begin{equation}
\label{mbc} Jf'+V f=\lambda f,\qqq
\begin{aligned}
&\qqq f_1(0)=\textstyle f_2(1)=0, \qq \{\t_n, n\in \Z\} \ mixed, \ 1 \ type
\\
& \qqq f_2(0)=\textstyle f_1(1)=0, \qq \{\vr_n, n\in \Z\} \ mixed, \ 2  \ type,
\end{aligned}
\end{equation}
where   $\lambda\in \C$. Here and in the following $f'$ denotes the derivative w.r.t. the first variable. The matrix $J$ and the matrix-valued potential $V$ are given by
\[
\lb{zse}
\begin{aligned}
 J= \ma 0&1\\-1&0\am,\qq
V= \ma v_1&v_2\\v_2&-v_1\am,\qq   \qq
v=\ma v_1 \\ v_2\am \in \mH.
\end{aligned}
\]
We assume that the vector $v$ belongs to the real Hilbert space $\mH=L^2((0,1),\R) \os L^2((0,1),\R)$,  equipped with the form
$\|v\|^2=\int_0^1 (v_1^2+v_2^2)dx$.
Let $\mu_n$ and $\nu_n, n\in \Z$ be
eigenvalues of the Dirichlet and Neumann problems respectively. Let
$\tau_n$ and $\varrho_n, n\in \Z$ be eigenvalues of the first and the
second problem  respectively with mixed boundary conditions
\er{mbc}, and we say shortly mixed eigenvalues.  All these
eigenvalues are simple and satisfy
\begin{equation}
\label{bx}
\begin{aligned}
& ......<\ol {\tau_1, \varrho_1}<\ol {\mu_1,\nu_1}<\ol {\tau_2, \varrho_2}< \ol
{\mu_2,\nu_2}<...,
\\
& \nu_n,\mu_n=\mu_n^o+o(1),\qq \tau_n, \varrho_n=\tau_n^o+o(1)\qq \as \ n\to
\pm \infty,
\end{aligned}
\end{equation}
where $\ol {u,v}$ denotes $\min \{u,v\}\le \max \{u,v\}$ for
shortness, and $\nu_n^o=\mu_n^o=\pi n$ and $
\tau_n^o=\varrho_n^o=\pi (n-{1\/2}), n\in \Z$ are the corresponding
unperturbed eigenvalues.
We introduce the fundamental solutions (vector-functions)  $\vt =(\vt_1, \vt_2)^\top$ and $\vp =(\vp_1, \vp_2)^\top$
of the equation $Jf'+Vf=\l f $,  under the conditions
$\vt (0,\l)=(1, 0)^\top$ and $ \vp (0,\l)=(0,1)^\top$.
Recall that $\m_n, \t_n$ and $\n_n, \vr_n, n\in \Z$ are zeros of the functions
$\vp_1(1,\l), \vp_2(1,\l)$ and $\vt_2(1,\l)$,  $\vt_1(1,\l)$ respectively.


We consider the operator $ T_{per} f=Jf'+Vf$ on $ L^2(0,2)\os L^2(0,2)$ with
2-periodic conditions $y(2)=y(0)$, where $v$ is 1-periodic and belongs to the real space $\mH$ on the unit interval.
The spectrum of $ T_{per} $ are eigenvalues $\l^\pm_{n}, n\in \Z$
which satisfy
$$
\begin{aligned}
& .....<\l^-_{1}\le \l^+_{1}<.... \le \l^+_{n-1}<
\l^-_n\le\l^+_{n}<...,
\\
& \l^\pm_{n}=\pi n+o(1)\qqq  \as \qq n\to \pm\infty.
\end{aligned}
$$
The eigenvalues $\tau_n, \varrho_n, $ and $\mu_n, \nu_n$ have the known  relations (see Fig. \ref{fig})
\begin{equation}
\label{exx}
\begin{aligned}
 \tau_n, \varrho_n\in (\l^+_{n-1},\l_{n}^-),\qqq  {\rm
and }\ \  \mu_n, \nu_n\in [\l^-_n,\l^+_{n}],\qqq \forall \ n\in \Z.
\end{aligned}
\end{equation}
  Here the equality $\l_n^-=\l_n^+$ means that $\l_n^-$ is a double
eigenvalue.  The eigenfunctions
corresponding to $\l_{n}^{\pm}$ have period 1 when $n$ is even and
they are antiperiodic, $y(x+1)=-y(x),\ x\in\R $, when $n$ is odd.
Recall that the operator $Jf'+Vf$ on the circle is the Lax operator for
 the periodic  defocusing Nonlinear Schr\"odinger equation (the NLS equation)
$iv_t=-v_{xx}+ 2|v|^2v$, see e.g., \cite{ZS72}, \cite{AT04}. The  NLS equation is one of
the most fundamental and the most universal nonlinear PDE.
Zakharov and Shabat proved that it   is integrable \cite{ZS72}.

\setlength{\unitlength}{1.0mm}
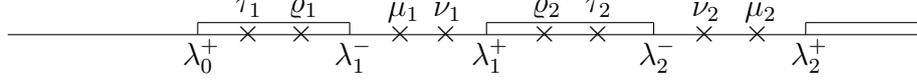
\begin{figure}[h]
\centering
\unitlength 1.0mm 
\begin{picture}(135,25)
\put(5,10){\line(1,0){120.00}}


\put(28,6){$\l_0^+$} \put(30,10){\line(0,1){1.6}}
\put(30,11.6){\line(1,0){20.00}} \put(50,10){\line(0,1){1.6}}
\put(48,6){$\l_1^-$}

\put(35,9){$\ts$} \put(42,9){$\ts$} \put(35,13){$\tau_1$}
\put(42,13){$\varrho_1$}

\put(55,9){$\ts$} \put(61,9){$\ts$} \put(55,12.5){$\mu_1$}
\put(61,12.5){$\nu_1$}

\put(66,6){$\l_1^+$} \put(68,10){\line(0,1){1.6}}
\put(68,11.6){\line(1,0){22.00}} \put(90,10){\line(0,1){1.6}}
\put(88,6){$\l_2^-$} \put(74,9){$\ts$} \put(81,9){$\ts$}
\put(74,13){$\varrho_2$} \put(81,13){$\tau_2$}

\put(95,9){$\ts$} \put(102,9){$\ts$} \put(95,12.5){$\nu_2$}
\put(102,12.5){$\mu_2$}

\put(108,6){$\l_2^+$} \put(110,10){\line(0,1){1.6}}
\put(110,11.6){\line(1,0){15.00}}

\end{picture}

\caption{\footnotesize \!\!\! Periodic $\l_n^\pm$, Dirichlet $\mu_n$,
Neumann $\nu_n$ and mixed $\tau_n$, $\varrho_n$ eigenvalues.} \label{fig}
\end{figure}

Introduce the real Banach spaces
$\ell^p=\ell^p(\Z), p\geq 1$ of real sequences $f=(f_n)_{n\in\Z}$ equipped with the norm $\|f\|_{(p)}^p=\sum |f_n|^p$.
Following the book of P\"oschel and Trubowitz \cite{PT87} we define
sets $\gJ^o, \gJ^1, \gJ$  of all real, strictly increasing sequences by
$$
\tes \gJ^o=\Big\{ s\!=\!(s_n)_{n\in \Z}: ....<s_1<s_2<....., \qq
s_n\!=\!\mu_n^o\!+\breve{s}_{n\,}   ,\qqq \breve{s}=(
\breve{s}_{n\,})_{n\in \Z}\!\in\!\ell^2     \Big\},
$$
$$
\tes
\gJ^1=\Big\{ s\!=\!(s_n)_{n\in \Z}: ....<s_1<s_2<....., \qq
s_n\!=\!\t_n^o\!+\breve{s}_{n\,},\qqq \breve{s}=(
\breve{s}_{n\,})_{n\in \Z}\!\in\!\ell^2     \Big\},
$$
$$
\tes \gJ=\Big\{ t\!=\!(t_n)_{n\in \Z}: ...<t_1<t_2<....., \qq
t_n\!=({\pi n\/2})^2+\breve{t}_{n\,},\qqq \breve{t}=(
\breve{t}_{n\,})_{n\in \Z}\!\in\!\ell^2 \Big\}.
$$
The mapping $s=(s_n)_{n\in \Z}\lra \breve{s} $ is a natural coordinate map between $\gJ^o$ and some
open convex subset $\breve{\gJ}^o=\Big\{ \breve{s}\!=\!(\breve{
s}_n)_{n\in \Z}\in\!\ell^2: ...<\mu_1^o+\breve{s}_{1\,}<\mu_2^o+\breve{s}_{2\,}<.....
\Big\} $ of $\ell^2$\,. Following  \cite{PT87} we identify $\gJ^o$
and $\breve{\gJ}^o$ using this mapping.  Below we refer to $\breve{s}=\!(\breve{
s}_n)_{n\in \Z}\in\!\ell^2$ as the standard coordinate system on $\cJ^o$.
As in \cite{PT87} this
identification allows to do analysis on $\gJ^o$ as if it was an open
convex subset of $\ell^2$. We have similar standard coordinate systems on $\cJ^1$ and $\cJ$.

Introduce  1-spectra mappings $\mu$ and $\nu$ from $\mH$ into $\gJ^o$
and $\t$ and $\varrho$ from $\mH$ into $\gJ^1$ by
\begin{equation}
\label{mntz1}
\begin{aligned}
v\to \mu=(\mu_n)_{n\in \Z},\ \qq  v\to \nu=(\nu_n)_{n\in \Z},
\qq
 v\to \t=(\tau_n)_{n\in \Z},\ \qq v\to \varrho=(\varrho_n)_{n\in \Z},
 \end{aligned}
\end{equation}
 For two 1-spectra mappings (only for strongly increasing and
alternate) we construct standard 2-spectra mappings of strongly
increasing sequences. For example, for $\t=(\tau_n)_{n\in \Z}\in \gJ^1$
and $\mu=(\mu_n)_{n\in \Z}
\in \gJ^o$ such that $....<\tau_1<\mu_1<\tau_2<\mu_2<...$
we define a 2-spectra mapping $\t \star \m$ from $\mH$ into $\gJ$ as
\begin{equation}
\label{st} v\to \tau \star \mu=(....,\tau_1,\mu_1,\tau_2,\mu_2,....).
\end{equation}

In order to describe inverse problems, following P\"oschel and Trubowitz \cite{PT87}, we introduce
the real norming constants $\gr_{n}, \gs_{n}$,  $\gt_{n}, \gu_{n}$
and the corresponding {\it norming } mappings by
$$
\begin{aligned}
& \gr_n=-\log |\vp_2(1,\m_n)|,\  \gs_n=-\log |\vt_1(1,\n_n)|,
\  \gt_n=-\log |\vp_1(1,\t_n)|, \
 \gu_n=-\log |\vt_2(1,\vr_n)|,
 \\
& v\to \gr =(\gr_n )_{n\in \Z},\ \  v\to \gs =(\gs_n )_{n\in \Z},\ \
v\to \gt  =(\gt_n )_{n\in \Z},\ \  v\to \gu   =(\gu_n )_{n\in \Z}.
\end{aligned}
$$
Recall that the vectors $(\m_n)_{n\in\Z}$ and $(\gr_n )_{n\in\Z}$ are
canonically conjugate variables for the NLS equation, see below Lemma \ref{T25} from  \cite{K01}. Note that it is  shown in \cite{FM76}
for the KdV equation.

We consider the four  spectra mapping, two spectra mapping, eigenvalues and norming constants mapping,  inverse periodic problems.
We describe  our main results:

\no $\bu$  We construct  the four  spectra mapping
and show that it is a real analytic bijection
between the space of potentials and the corresponding spectral data.

\no $\bu$  The basic inverse problems are isomorphic and the corresponding automorphisms are obtained in explicit forms.
Each of these inverse problems   is a real analytic bijection
between the space of potentials and the corresponding spectral data.

\no $\bu$ We define and describe new inverse problems: shifting, replacing mappings.


To the best of our knowledge the obtained results have no analogies
in existing literature. We need to underline that in order to discuss isomorphic inverse problems for ZS-systems on the unit interval we need also
results about inverse periodic problems.

Our proof uses observations 1)-3) and also
following results and methods about ZS-systems:


\no $\bu$ Uniqueness for inverse problems, see e.g. \cite{A14}.

\no $\bu$ The spectral parameters $\n, \gs$ are locally free, see \cite{DK00}
(the explicit transforms corresponding to the change of only a finite number of spectral parameters, eigenvalues plus norming constants).

\no $\bu$ The mapping $\m\ts \gr$ is a real analytic local isomorphism, see e.g., \cite{K01}, \cite{K05}.

\no $\bu$ Inverse periodic problem (a characterization, a priori estimates)
from \cite{K01}, \cite{K05}, \cite{K96}.

\subsection{Main results on the unit interval}
In order  to discuss main results we need new basic inverse problem
via four 1-spectra mappings ZS-operators.
Following recent
paper \cite{K19} we define the 4-spectra mapping $\gf: \mH\to \ell^2$  by
\begin{equation}
\label{dgf}
\tes
v\to \gf(v)=(\gf_n(v))_{n\in \Z},\qq \gf_{2n-1}={1\/2}(\varrho_n
-\tau_n),\qqq  \gf_{2n}={1\/2}(\nu_n -\mu_n), \ \ \ \ \ n\in \Z.
\end{equation}
We sometimes write $\mu_n(v), \nu_n(v),...$
instead of $\mu_n, \nu_n,...$,  when several potentials  are being
dealt with.  Recall some definitions. We write $\mathcal K_\C$ for the complexification of the real Hilbert space $\mathcal K$. Suppose that $ \mathcal K, \mathcal S$ are real separable
Hilbert spaces. The mapping $f:\mathcal K\to \mathcal S$ is a local real analytic isomorphism iff for any $y\!\in\!\mathcal K$ it has an analytic
continuation $\widetilde {f}$ into some complex neighborhood
$\!\V\!\ss\!{\mathcal K}_\C$ of $ y$, which is a bijection between $\V$ and
some open set $\widetilde {f}(\V)\!\ss\!{\cS}_\C$ and if $\widetilde {f}$,
$\widetilde {f}^{-1}$ are analytic mappings on $\V$\,, $\widetilde {f}(\V)$
respectively. The mapping $f$ is a real-analytic bijection (shortly
a RAB) between $\mathcal K$ and $\mathcal S$ if it is both a bijection and a local real analytic isomorphism.

\begin{theorem}
\label{T1}
The 4-spectra mapping $\gf: \mH\to \ell^2$ defined by \er{dgf},
 is a RAB between $\mH$ and $\ell^2(\Z)$ and  satisfies
\begin{equation}
\label{egf}
\begin{aligned}
 \tes
 {1\/ \sqrt 2}\|\gf(v)\|\leq \|v\|\leq 2\|\gf(v)\|\big(1+\|\gf(v)\|\big)\qqq
 \forall \ v\in\mH.
\end{aligned}
\end{equation}

\end{theorem}

\no {\bf Remark.} Define the set $\mD$ of all $v\in\mH$ such the sequence
$(\gf_n(v))_{n\in\Z}$ is finitely supported, i.e.,  $\gf_n(v)=0$ for all $n\in\Z$ large enough. The  4-spectra mapping $\gf: \mH\to \ell^2$  is a RAB between $\mH$ and $\ell^2(\Z)$ and then the set $\mD$
 is dense in $\mH$.

 Let $\gS$ be a set of all diagonal operators $\sigma= \diag
(...,\s_1,\s_2,...)$ on $\ell^2$, or shortly $\s=(\s_j)_{n\in \Z}$, where
$\s_j\in \{\pm 1\}, j\in \Z$. This set $\gS$ defines the so-called
lamplighter group, see \cite{CM17}. For each $\s\in \gS$ and $\gf$ is given by \er{dgf}, we define a lamplighter mapping $ \mathcal U_\s$ by
\begin{equation}
\label{dj}
 \cU_\s=\gf^{-1} \ci(\s \gf): \mH\to
\mH.
\end{equation}
Note that $\cU_\s:\mD\to \mD$ for all
$\s\in\gS$.
We define a reflection  $\cR$  and reflection type operators $\cF_0,\cF_1, \cF_2$ acting on $\mH$ by
\[
\lb{rA}
(\cR v)(x)= v(1-x),x\in[0,1],\qq {\rm and}\qq
\cF_0=-\1_2, \ \ \cF_1=J_1\cR,  \qqq      \cF_2=-\cF_1.
\]

\begin{theorem}
\label{T2}

i) Each mapping $\cU_\s: \mH\to \mH,  \s\in \gS$  is a RAB of $\mH$ onto itself and  satisfies
\begin{equation}
\label{dUx}
\begin{aligned}
 \mathcal U_\s=\mathcal U_\s^{-1}, \qqq \mathcal U_\s\ci \mathcal U_{\s'}=\mathcal U_{\s\s'}=\mathcal U_{\s'}\ci
\mathcal U_{\s},\qqq \forall \qq \s, \s'\in \gS,
\end{aligned}
\end{equation}
\begin{equation}
\label{dUq} \|\cU_\s(v)\|=\|v\|\qqq \forall \ v\in\mH.
\end{equation}
\no ii) The 2-periodic   eigenvalues $(\l_{n}^\pm)_{n\in\Z}$
are invariant under each $\cU_\s, \s\in \gS$, i.e.,
\[
\label{peU} (\l_{n}^\pm)_{n\in\Z}=(\l_{n}^\pm)_{n\in\Z}\ci \cU_\s
\qq \forall \ \s\in \gS,
\]
and the mappings $\cU_\s$ for  specific $\s\in\gS$ have the forms:
\begin{equation}
\label{dU}
\begin{aligned}
\cF_o=\mathcal U_\s,\qq \where \ \s=-I,\qqq
\ca
\cF_1=\cU_\s,\qq \where \qq \s_{n}= (-1)^n\qqq
\forall \ n\in \Z
\\
\cF_2=\cU_\s,\qq \where \qq \s_{n}= -(-1)^n\qq
\forall \ n\in \Z\ac.
\end{aligned}
\end{equation}

\no iii)   If  $\s=(\s_j)_{n\in \Z}\in \gS$, then for each
 $n,j\in \Z$ we have
 \begin{equation}
\label{U1}
{\rm if} \ n=2j-1 \Rightarrow \ca (\tau_j, \vr_j)=(\varrho_j,\t_j)\ci \cU_\s, \ \ {\rm if}\qq  \s_n=-1
\\
(\tau_j, \vr_j)=(\t_j,\varrho_j)\ci \cU_\s, \ \ {\rm if} \qq
\s_n=1 \ac,
\end{equation}
\begin{equation}
\label{U2} {\rm if} \ n=2j \Rightarrow \ca (\m_j, \n_j)=(\n_j,\m_j)\ci \cU_\s, \ \ {\rm if}\qq  \s_n=-1
\\
              (\m_j, \n_j)=(\m_j,\m_j)\ci \cU_\s
              \ \ {\rm if} \qq \s_n=1 \ac .
\end{equation}

\end{theorem}

\no {\bf Remark.} 1) The mapping $\gf$ is non-linear, but $\mathcal U_\s$
keeps the norm on $\mH$, see \er{dUq}.

\no 2) Due to \er{dU} the mapping $\cU_\s$ has the very simple form
for specific operators $\s$.

We discuss inverse problems for 1-spectra mappings and  norming mappings.

\begin{theorem}
\label{T3}

\no i) All 2-spectra mappings $\t\star \m, \varrho\star \n$,  $\varrho\star
\m$ and $\t\star \n$ acting from $\mH$ into $\gJ$  are isomorphic, each of them  is a RAB between $\mH$ and $\gJ$ and they satisfy
\begin{equation}
\label{enm}
\begin{aligned}
\t\star \m=(\varrho\star \n)\circ \cF_o=(\varrho\star \m)\circ \cF_1=(\t\star
\n)\circ \cF_2.
\end{aligned}
\end{equation}
\no ii)  Each of mappings $\m\ts \gr , \n \ts \gs $,  $\m \ts \gs $
and $\n \ts \gr $ acting from $\mH$ into $\gJ^o\ts \ell^2$
is a RAB between $\mH$ and $\gJ^o\times\ell^2(\Z)$
Moreover, they are isomorphic and satisfy
\[
\label{enm1} \begin{aligned}
\m\ts \gr =(\n \ts \gs )\circ \cF_o
=(\m \ts (-\gr) )\circ \cF_1=(\n \ts (-\gs) )\circ \cF_2.
\end{aligned}
\]
\no iii)  Each of mappings $\t\ts (\pm\gt), \vr \ts (\pm\gu ) $,
acting from $\mH$ into $\gJ^1\ts \ell^2$
is a RAB between $\mH$ and $\gJ^1\times\ell^2(\Z)$.
Moreover, they are isomorphic and satisfy
\[
\label{enm2} \begin{aligned}
\t\ts \gt =(\vr \ts \gu )\circ \cF_o
=(\vr \ts (-\gu) )\circ \cF_1=(\t \ts (-\gt) )\circ \cF_2.
\end{aligned}
\]
\end{theorem}

This theorem shows that the mapping $\m\ts \gr , \n \ts \gs $,... are isomorphic. But in order to show that the mappings
$(\cS\m)\ts \gr, \t \ts \gt$ are isomorphic we need to introduce a shifting mapping.
We discuss new inverse promblems about a {\it shifting mapping}
$\cS:\cJ^o\to\cJ^1$ defined by
\[
\lb{A3}
\tes
 (\cS z)_n=z_n-{\pi\/2},\qq z=(z_n)_{n\in \Z}\in \cJ^o, \  n\in \Z.
\]

\no {\bf Problem}: {\it
Consider the Dirichlet eigenvalues $\m_n(v), n\in\Z$ for some $v\in\mH$
and a sequence $\cS\m(v)$. Do we have $u\in \mH$ such that $\t(u)=\cS\m(v)$? Can we describe $u\in \mH$?}

In order to study such problem
we define the mapping $\gF:\mH\to \mH$ by $\gF v(x)=e^{{\pi}x J}v(x)$.

 \begin{theorem}\label{T4} 
 Each of the mappings
$(\cS\m)\ts \gr, (\cS\n)\ts \gs, \t \ts \gt$ and $\vr \ts \gu$ acting from $\mH$ into $\gJ^1\ts \ell^2$ is a RAB between $\mH$ and $\gJ^1\ts \ell^2$  and they satisfy
\begin{equation}
\label{enz1}  (\cS\m)\ts \gr=(\t \ts \gt)\circ \gF=((\cS\n)\ts \gs)\circ \cF_o=(\vr \ts \gu)\ci \cF_o\circ \gF,
\end{equation}
where  the mapping $\cS:\cJ^o\to \cJ^1$ is a bijection between
$\cJ^o$ and $\cJ^1$.

 \end{theorem}
 \no {\bf Remark.}  By this theorem, the  mappings  $\m\ts \gr, \t\ts \gt$ are isomorphic and we do not know such effect for the case of Schr\"odinger operators on the unit interval, see \cite{K25}.

Define the even-odd space $\mH_{eo}$ and the odd-even space $\mH_{oe}$ by
$$
\mH_{eo}=\big\{v\in \mH: v=J_1\cR v\big\}, \qqq
\mH_{oe}=\big\{v=(v_1,v_2)^\top: (v_2,v_1)^\top\in \mH_{eo}\big\}.
$$

  \begin{corollary}
\label{T5}
Each of the mappings $\cS\m, \cS\n, \t$ and $\vr$ acting from $\mH_{eo}$ into $\gJ^1$  is a RAB between $\mH_{eo}$ and $\gJ^1$ and they satisfy on the space $\mH_{eo}$:
 \begin{equation}
\label{enz3}
\cS\m=(\cS\n)\ci \cF_o=\t\circ \gF=\vr\ci \cF_o\circ \gF,
\end{equation}
where  the mapping $\gF :\mH_{eo}\to \mH_{oe}$ is a bijection
between $\mH_{eo}$ and $\mH_{oe}$.
 \end{corollary}

 By this theorem, the  mapping  $\m, \n, \t$ and $\vr$ are isomorphic on the space $\mH_{eo}$ and we do not know such effect for the case of the Schr\"odinger operators on the unit interval $[0,1]$, see \cite{K25}.

{\bf Replacing mappings.}
 We discuss a new type of inverse problems.  Let the Dirichlet
mapping $v\to \m=(\mu_n)_{n\in \Z}$ be given  and replace some $\mu_n$ by
the Neumann eigenvalues $\nu_n$. Then we obtain a replacing mapping
$\gc$. For example, we have $\gc=(...,\mu_1,\nu_2, \nu_3,\mu_4, \mu_5,...)$.
There is a question: it is a good 1-spectra mapping?
We discuss {\it replacing} mappings on the finite intervals.
Let $\Y_1, \Y_2$  be some subsets of $\Z$.
We define replacing mappings $ v\to\z=(\z_n)_{n\in\Z},
v\to\phi=(\phi_n)_{n\in\Z}$ and
their components by
\[
\lb{Rep}
\begin{aligned}
\z_n=\left\{\begin{aligned}
&\vr_n, \qq \text{if}  &&n\in\Y_1 \\
&\t_n, \qq \text{if}  &&n\notin\Y_1
\end{aligned}\right., \qq
\qq
\phi_n=\left\{\begin{aligned}
&\n_n, \qq \text{if}  &&n\in\Y_2 \\
&\m_n, \qq \text{if}  &&n\notin\Y_2
\end{aligned}\right..
\end{aligned}
\]
Define the operator $\s=(\s_n)_{n\in\Z}$ by
\[
\lb{Rep1}
\s_{2j-1}=\ca  -1   &  \text{if} \ j\in \Y_1 \\
               1  & \text{if} \ j\notin \Y_1   \ac,\qq
               \
               \
               \s_{2j}=\ca  -1   &  \text{if} \ j\in \Y_2 \\
               1 & \text{if} \ j\notin \Y_2   \ac.
\]

\begin{corollary} \lb{Ter1}
Define  the mappings $\z=(\z_n)_{n\in\Z}$ and
$\phi=(\phi_n)_{n\in\Z}$ and $\gz=(\gz_n)_{n\in\Z}$  by \er{Rep}.  Then
the mapping $\z\star \phi: \mH\to\gJ $ is a RAB between $\mH$ and $\gJ$ and satisfies
\begin{equation}
\label{rtm}
(\z\star \phi)\circ \cU_\s=\t\star \m.
\end{equation}

\end{corollary}

\no {\bf Final remarks.} 1)  
 Theorems \ref{T3}, \ref{T4} and asymptotics of $\m_n, \gr_n$ give
 the asymptotics of eigenvalues $\m_n,..,\r_n$ and the norming constants $\gr_n,..,\gu_n$ as $n\to \pm\iy$, see Theorem \ref{Tas}.
  
 \no 2) The canonical relations between $\m_n, \gr_n, n\in \Z$ from Lemma \ref{T25} and isomorphisms of inverse problems give that the corresponding pairs $(\n_n, \gs_n), (\t_n, \gt_n)$ and $(\vr_n, \gu_n), n\in \Z$ are also canonical  variables, see Theorem \ref{Tcv}.

 \no 3) The case of periodic ZS systems is discussed in
Section 5.


\subsection{Short review}
We shortly  describe of known results in the inverse spectral theory for differential operators on the  unit interval.
These inverse problems   were investigated by many authors
(Borg, Gel'fand, Levitan, Marchenko, Trubowitz, ..), see  the monographs \cite{L87}, \cite{M86}, \cite{PT87} and references therein.
 We recall only some important steps mostly focusing on the {\it characterization} problem.
Borg obtained the first result about uniqueness for two spectra mapping,
improved by \cite{Le49}.
Marchenko \cite{M50} proved that the Dirichlet eigenvalues
plus the  normalizing constants determine the
the potential uniquely. Gel'fand and Levitan \cite{GL51} created a basic method to reconstruct this potential via the famous integral equation.
Remark that independantly, a different approach to this problem was developed by Krein \cite{Kr51}, \cite{Kr54}.
 At that time, there was a gap between necessary and sufficient conditions
for inverse problems corresponding to fixed classes of potentials.
 Marchenko and Ostrovski \cite{MO75} gave the complete
solution of the inverse problem in terms of two spectra. Note that some results about that were obtained by Levitan and Gasymov \cite{LG64}.
Trubowitz and co-authors, \cite{IT83}, \cite{IMT84}, \cite{PT87}, suggested another approach.
It is based on analytic properties of the mapping
$\mathrm{\{potentials\}\mapsto\{spectral\ data\}}$ and the explicit transforms corresponding to the change of only  a {\it finite} number of spectral parameters (eigenvalues plus norming constants). This approach was developed in \cite{CM93},  \cite{CM93a}, \cite{KC09} and
was applied to other inverse problems with purely discrete spectrum:
(a) for periodic case, see \cite{GT87}, \cite{K97}, \cite{K99}, (b)
perturbed harmonic oscillator, see  \cite{MT81}, \cite{CKK04}, \cite{CK07},\
(c) Sturm-Liouville problems with matrix-valued potentials under the Dirichlet  boundary conditions on the unit interval, see  \cite{CK06}, \cite{CK09}.

The inverse problems for ZS-operators on a finite interval also are well  studied.  Uniqueness and the reconstruction results were obtained in \cite{GL66}, \cite{LS91}, see also \cite{Ge02}, \cite{S70}.  Explicit formulas for solutions (based on the degenerate Gelfand–Levitan equation) in the case where finitely many spectral data are perturbed were given in \cite{DK00}. Uniqueness results for other types of inverse problems were established, see e.g., for mixed spectral \cite{H01} or interior  data \cite{DG01}, \cite{MT02}. Misura  solved inverse problems for 2-spectra mapping  $\t\star \m$ (including the characterization) in \cite{Mi78}, \cite{Mi79}. The proof  is essentially the same as in  \cite{MO75}. Later on the characterization problems  (for  inverse problems for a 2-spectra mapping  $\vr\star \n$ and the mapping $\gu\star \vr$ was solved in \cite{AHM05}.

There are a lot of results about the periodic case, see e.g., \cite{Mi78}, \cite{Mi79}, \cite{Mi80}, \cite{GG93}, \cite{GK14}, \cite{K01},
\cite{K05}, \cite{K01}, \cite{K96} and the references therein.
Misura  solved inverse problem in terms of conformal mapping (including the characterization) in \cite{Mi78}, \cite{Mi79}. The proof  is essentially the same as in  \cite{MO75}. Korotyaev \cite{K01} solved the inverse problems
(including characterization and stability estimates) in term
 of the local maxima and minima of Lyapunov functions on the real line.
In the next papar \cite{K05} he solved the
inverse problems  in terms of gap lengths. The proof (including
characterization) was based on the analytic method from \cite{GT87}, \cite{K97}. In this approach a
priori (or stability) estimates of potentials in terms of spectral data are crucial.

We discuss the stability estimates. We consider only sharp cases, which are obtained only for periodic case.
The two-sided estimates of a potential in terms of  gap lengths (or parameters
of the Lyapunov function) were obtained in \cite{K96}, \cite{K98}, \cite{K00}, \cite{K05}, \cite{K06}
via the conformal mapping theory  associated with quasimomentum.
Here results about various properties of the conformal mapping theory
from \cite{MO75}, \cite{KK95}, \cite{K97m}.
 Recall that Hilbert \cite{H09} obtained the first result about such  conformal mappings from a multiply
connected domains onto a domain with parallel slits.


\section {Preliminary results \lb{Sec2}}
\setcounter{equation}{0}


\subsection{Fundamental solutions}
 Let $J_1= \ma 1&0\\0&-1\am,\ \         J_2= \ma 0&1\\1&0\am . $
Below we need the simple identities  for $J, J_1, J_2, V$ and all $\l\in \C$:
\[
\lb{11}
  J^2=-I, \ \ \ J_1J_2=J,\ \ \ JJ_1=-J_2,\ \ \ JJ_2=J_1,
\]
\[
\lb{12}
e^{\l J}=\1_2\cos \l +J\sin \l,
\]
\[
\lb{13}
JV=-VJ,\qqq e^{\l J} V=Ve^{-\l J}.
\]
For each $(\l,v)\in \C\ts \mH_\C$ the $2\times 2$- matrix valued solution of the following equation
\[
\lb{15}
Jy'+Vy=\l y , \ \ y(0,\l )=\1_2,\ \ x\ge 0,
\]
 has the representation
\[
\lb{16}
y (x,\l)=e^{-\l xJ}+\int_0^xe^{-\l J(x-t)}JV(t)y(t,\l)dt.
\]
It is known that the solution of this integral equation has the following form
\[\lb{17}
y(x,\l)= \sum _{n\geq 0} y_n(x,\l),\qqq y_0(x,\l)=e^{-\l xJ}=\ma \cos \l x&-\sin \l x\\ \sin \l x&\cos \l x\am,
\]
where the functions $y_n$ are defined by the relations:
\[\lb{18}
y_n(x,\l)=\int_0^xe^{-\l J(x-t)}JV(t)y_{n-1}(t,\l)dt,\ \ n\geq 1.
\]
For fixed $x,\l$, the function $y_n $ is a multi-linear form on
$\mH_\C\ts\dots \ts \mH_\C$. Using \er{18} we have
$$
y_1(x,\l,v)=\int_0^xe^{-\l J(x-t)}JV(t)e^{-\l tJ}dt=\int_0^xe^{-\l J(x-2t)}JV(t)dt,
$$
and so on.  We have the Wronskian identity
$$
\det y =\vt_1\vp_2-  \vt_2\vp_1=1.
$$
In order to present results from \cite{K05} we define the functions
$$
(a,b)_j:= (J_ja,b),\ \ \ a,b\in \R^2 , \ \  j=0,1,2,\qq J_0:=J.
$$

\begin{lemma}  \lb{T21}
i) For each $(\l,v)\in \C\ts \mH_\C$
there exists  a unique solution $y$ of Eq. \er{16} having the form
\er{17}-\er{18}, where series \er{17} converge
absolutely and uniformly on bounded subsets of $[0,1]\times\C\times \mH_\C$.
For each $x\in [0,1]$ the function
$y(x,\l,v)$ is entire on $\C\times \mH_\C$ and satisfies
\[
 \lb{fuso}
|y(x,\l,v)|\leq e^{ |\Im \l|x+\|v\|_\C}.
\]
Moreover, if the sequence $v^{(s)}$ converges weakly to $v$ in $\mH_\C$, as $s\to \iy$, then $y (x,\l,v^{(s)})\to y (x,\l,v)$ uniformly
on bounded subsets of $[0,1]\times \C$.

\no ii) The derivatives of $y(1,\l,v)$ with respect to $v=(v_1,v_2)^\top$ have the forms
$$
{\pa y(1,\l,v)\/\pa v_j(x)}=
\ma
\wt\vt_1(\vt,\vp )_j-\wt\vp_1(\vt,\vt )_j&
\wt\vt_1(\vp,\vp )_j-\wt\vp_1(\vt,\vp )_j\\
\wt\vt_2(\vt,\vp )_j-\wt\vp_2(\vt,\vt )_j&
\wt\vt_2(\vp,\vp )_j-\wt\vp_2(\vt,\vp )_j
\am ,\qq j=1,2.
$$
where $\vt =\vt (x,\l,v), \vp =\vp (x,\l,v)$,
 and $\wt\vt =\vt (1,\l,v), \wt\vp =\vp (1,\l,v)$.

\no iii) The following asymptotics hold true as $|\l|\to \pm\iy $:
\[ \lb{2.19}
y(x,\l,v)=y_o(x,\l)+ o(1)e^{|\Im \l|x} ,
\]
 uniformly on the bounded subsets of $[0,1]\times\C \ts
\mH_\C$.
\end{lemma}


\subsection{Dirichlet eigenvalues and norming constants}
We recall needed results about
Dirichlet eigenvalues $\m_n$ and the corresponding norming constants
$\gr_n , n\in \Z$ from \cite{K01}, \cite{K05}.

\begin{lemma}  \lb{T23}

i) Let  $v\in\mH_\C$ and $\ve_v =4^{-4}e^{-3\|v\|}$ and $m_v\in \N$ be large enough. Then for each
integer $m>m_v$ and any $u\in \cB_\C(v,\ve_v) $ the function $\vp _1(1,\l,u)$
 has exactly $2m+1$ roots, counted with multiplicities,
in the disc $\{z: |z|<\pi (m+\tes{1\/2})\}$
and for each $|n|>m,$ exactly one simple root in the disc
$\{z: |z-\pi n|<1\}.$  There are no other roots.

\no ii) The 1-spectra mapping $v\to \m=(\m_n(v))_{n\in \Z}$ is a real analytic from $\mH$ in $\gJ^o$. Furthermore, each function $\m_n, n\in \Z,$ is compact and  real analytic  on $\mH$ and its gradient is given
\[
{\pa \m_n(v)\/\pa v(x)}=
{\big((\vp , \vp )_1,(\vp , \vp )_2\big)
 \over \|\vp(\cdot,\m_n(v),v)\|^2} (x,\m_n(v),v).
\]
iii) The mapping $v\to \gr=(\gr_n (v))_{n\in \Z}$ is a real analytic from $\mH$ into $\ell^2$. Furthermore, each function   $\gr_n=-\log |\vp_2(1,\m_n,\cdot)|, n\in \Z$, is compact, real
analytic on $\mH$ and its gradient is given  by
\[
{\pa \gr_n (v)\/\pa v(x)}=-(-1)^n e^{\gr_n (v)}
\Big(\dot \vp_2(1,\l,v){\pa \m_n(v)\/\pa v(x)}+{\pa \vp_2(1,\l,v)\/\pa v(x)}\Big)\Big|_{\l=\m_n(v)}.
\]

\end{lemma}

Recall that $(a,b)_0=a_1b_2-a_2b_1$, for $a,b\in\C^2$.
Define the symplectic form
$$
f\wedge g=\int _0^1 (f(x),g(x))_0dx,\ \ \ f,g\in \mH,
$$
and note that $f\wedge f=0$.
Below we need  the following canonical relations from \cite{K01}.

\begin{lemma}  \lb{T25}
 For any $n, j\in \Z$, the following identities for $\m_n, \gr_n, n\in\Z$ hold true:
\[
\lb{mr}
\m_n'(v)\we \m_j'(v)=0,\qqq
\gr'_{n}(v)\we \m_j'(v)=\d_{nj},
\qqq
\gr'_{n}(v)\we \gr'_{j}(v)=0,
\]
where $\m_n'={\pa\m_n\/\pa v},\gr'_{n}={\pa \gr_n \/\pa v}$ and the sequence $\{ \m_n', \gr'_{n}, \ n\in \Z \}$,
is a basis for $\mH$.
\end{lemma}

Introduce the Fourier transformation
$\F : \mH_C\to \ell^2_C\os \ell^2_C \ $ by the formulas   $\F f=((\F f)_n)_{n\in \Z}$, where the components $\hat f_n:=(\F f)_n$ have the forms
$$
 \hat f_n= \ma f_{1n}\\  f_{2n} \am =\int_0^1
e^{2\pi nxJ}f(x)dx= \ma f_{1(nc)}+f_{2(ns)}\\ - f_{1(ns)}+f_{2(nc)} \am,
\ \ f(x)=\ma f_{1}\\  f_{2} \am (x)=\sum_{n\in \Z} e^{-2\pi nxJ}\hat f_n,
$$
$$
f_{j(nc)}=\int_0^1f_j(x)\cos 2\pi nx dx, \qqq
f_{j(ns)}=\int_0^1f_j(x)\sin 2\pi nx dx, \qqq  j=1,2.
$$
In analogy to the notation $O(1/n)$ we use the notation
$\ell^d(n),d\geq 1$, for an arbitrary sequence of numbers, which is
an element of $\ell^d$ (see \cite{PT87}). For instance,
$$
 a_n=b_n+ \ell^d(n)\ \ \ {\rm is \ equivalent  \ to}
\ \ \ (a_n-b_n)_{n\in\Z}\in \ell^d .
$$
We recall results  from \cite{K01}, \cite{K05} about the  mapping
$\m\ts \gr$, where $\m_n$ are the Dirichlet eigenvalues and $\gr_n, n\in \Z$ are their norming constants.

\begin{theorem}
\label{Tdnc}
The mapping $\m\ts \gr: \mH\to \cJ^o\ts \ell^2$
is a real analytic local isomorphism between $\mH$ and $\cJ^o\ts \ell^2$.
Moreover, for any fixed $d>1$ the following asymptotic
estimates hold true:
\[
\lb{ra1}
F_n(v):=\ma \m_n (v)\\ \gr_n (v)\am=\ma \pi n\\ 0\am -J_1(\F v)_n+\ell^d(n),
\]
\[
\lb{ra2}
{\pa F_n(v)\/\pa v(x)}=-J_1  {\pa (\F v)_n\/\pa v(x)}+\ell^ 2(n),
\]
 as $n\to\pm \infty $, uniformly on
$[0,1]\times \cB_\C(u,\ve_u)$ for each $u\in \mH$, where $\ve_u=
4^{-4}e^{-3\|u\|}$.

\end{theorem}

\no {\bf Proof.} By Lemma \ref{T23},
the mapping  $F=\m\ts \gr: \mH\to \cJ^o\ts \ell^2$
is real analytic, and the asymptotics \er{ra1}, \er{ra2} hold true.
We show that this mapping $F$ is a local real analytic isomorphism.
Let $F'(v)={\pa F(v)\/\pa v}$ for shortness.
By \er{ra2}, for fixed $v\in\mH$ the operator  $ F'(v)- F'(0)$ is a compact
and $ F'(0)$ is the Fourier transformation from $\mH$ onto $\ell^2\os \ell^2$. Thus $F'(v)$ is a Fredholm operator.
We prove that the operator $F'(v)$ is invertible.

Let $\x\in \mH, \x\ne 0$ be a solution of the equation
$$
F'(v)\x=0\ \ \ \  \Leftrightarrow \ \ \ \ \   \ \lan \m_n'(v),\x\ran=0, \ \   \lan \gr_n (v),\x\ran=0    \qq \forall \ n\in \Z.
$$
Due to Lemma \ref{T25}, the sequence $\{ \m_n'(v), \gr_n'(v), \ n\in \Z \}$
is a basis for $\mH$. Then we obtain $\x=0$, which gives the contradiction
and the operator $F'(v)$ is invertible. Thus
$F$ is the local real analytic isomophism.
\BBox

The functions $\textstyle \vp_j(1,\l), \vt_j(1,\l), j=1,2$ are entire and have the Hadamard factorizations
\[
\lb{FS}
\begin{aligned}
&\vp_1(1,\l)=(\m_0^o-\l) \ {\rm v.p.}{\prod_{n\in \Z,n\ne0}} {\mu_n-\l\/\mu_n^o}, \qq \vp_2(1,\l)=  {\rm v.p.}\prod_{n\in \Z}{\t_n-\l\/\t_n^o},
\\
&\vt_1(1,\l)={\rm v.p.}\prod_{n\in \Z}{\vr_n-\l\/\vr_n^o}, \qq \vt_2(1,\l)=(\n_0-\n_0^o) \ {\rm v.p.}{\prod_{n\in \Z,n\ne0}} {\nu_n-\l\/\nu_n^o},
\end{aligned}
\]
where the products converge uniformly on compact sets on the complex plane.
Following Marchenko \cite{M50}, we define
the real normalizing constants $\ga_n , \gb_n $,  $\gc_n , \gd_n , n\in \Z$
for correponding eigenvalues $\m_n, \n_n, \t_n, \vr_n$ respectively and the corresponding {\it normalizing } mappings
by
\begin{equation}
\label{mnc}
\begin{aligned}
& \tes e^{-\ga_n}= \|\vp(\cdot,\mu_n)\|^2,
\
e^{-\gb_n}=\|\vt(\cdot,\nu_n)\|^2,
\
e^{-\gc_n}=\|\vp(\cdot,\t_n)\|^2,\  e^{-\gd_n}= \|\vt(\cdot,\vr_n)\|^2,
\\
& v\to \ga =(\ga_n )_{n\in \Z},\   v\to \gb =(\gb_n )_{n\in \Z},\
v\to \gc  =(\gc_n )_{n\in \Z},\   v\to \gd   =(\gd_n )_{n\in \Z}.
\end{aligned}
\end{equation}
Note that  $\ga(0)=\gb(0)=\gc(0)=\gd(0)=0$.
We rewrite the
 the normalizing constants $\ga_n , \gb_n $,  $\gc_n , \gd_n$ in terms of  the norming constants $\gr_n , \gs_n $,  $\gt_n , \gu_n , n\in \Z$
in the following forms
\begin{equation}
\label{ss1}
\begin{aligned}
\tes
& e^{-\ga_n }=-\dot \vp _1(1, \m_n)\vp_2(1,\m_n)
 =e^{ \ga_n ^\bu -\gr_n },\qq  \ga_n ^\bu=\ln |\dot \vp _1(1, \m_n)|,
\\
&  e^{-\gb_n}=\dot \vt_2(1,\n_n)\vt_1(1,\n_n)
=e^{\gb_n^\bu -\gs_n },\qqq \gb_n^\bu=\ln |\dot \vt_2(1,\n_n)|,
\end{aligned}
\end{equation}
\begin{equation}
\label{ss2}
\begin{aligned}
e^{-\gc_n}=-\dot \vp_2(1, \t_n)\vp_1(1,\t_n)=e^{\gc_n^\bu-\gt_n },
\qq \gc_n^\bu=\ln |\dot \vp_2(1, \t_n)|,
\\
e^{-\gd_n}=\dot \vt_1(1,\vr_n)\vt_2(1,\vr_n)
 =e^{  \gd_n^\bu   -\gu_n },\qq \gd_n^\bu=\ln |\dot \vt_1(1,\vr_n)|,
\end{aligned}
\end{equation}
where $\dot u={\pa u\/\pa \l}$. In order to discuss $\ga_n ^\bu,...,\gd_n ^\bu$ we recall that for all $n\in\Z$ we have:
$$
(-1)^n\vp_2(1,\m_n)>0,\ \ (-1)^n\vp_1(1,\n_n)>0,\ \ (-1)^n\vp_1(1,\t_n)>0, \ \ (-1)^{n+1}\vt_2(1,\vr_n)>0.
$$

\begin{proposition} \label{Tck09}
Let $v\in\mH$. Then following asymptotics  hold true:
\begin{equation}
\label{aa1}
\begin{aligned}
\ga_n^\bu(v)=\ell^2(n),\qq \gb_n^\bu(v)=\ell^2(n),\qq \gc_n^\bu(v)=\ell^2(n),\qq \gd_n^\bu(v)=\ell^2(n),
\end{aligned}
\end{equation}
as $n\to \pm\infty$ uniformly on $\cB_\C(u,\ve_u)$ for each $u\in \mH$, where $\ve_u=4^{-4}e^{-3\|u\|}$.

\end{proposition}

\no{\bf Proof.} We have $\vp_2^o(1,\l)=\cos \l$ and
$\dot\vp_2^o(1,\m_n^o)=(-1)^n$ at $v=0$.
From \er{ra1} we have $\t_n=\t_n^o+\breve{
\t}_n$, where $(\breve{\t}_n)_{n\in \Z}\in \ell^2$ and let $a_{j,n}=\frac{\breve{\t}_j-\breve{\t}_n}{\t_j^o-\t_n^o}$.
Then using \er{FS} we obtain
$$
(-1)^n\dot \vp_2(1,\t_n)={\rm v.p.}\prod_{ j\ne
n}\frac{\t_j-\t_n}{\t_j^o}
={\rm v.p.}\prod_{j\ne
n}\frac{\t_j\!-\!\t_n}{\t_j^o-\t_n^o}
={\rm v.p.}\prod_{j\ne n}\Big[1+a_{j,n}\Big],\ \
$$
which yileds
$$
\label{fS}
\begin{aligned}
& \log[(-1)^n \dot \vp_2(1,\t_n)]=\log{\rm v.p.}\prod_{j\ne
n}\Big[1+a_{j,n}\Big] =
{\rm v.p.}\sum_{j\ne n}\Big[a_{j,n}+O(a_{j,n}^2) \Big]
\\
& ={\rm v.p.}\sum_{j\ne n}\frac{\breve{
\t}_j-\breve{\t}_n}{\pi(j-n)}+\ell^1(n)=
{\rm v.p.}\sum_{j\ne n}\frac{\breve{
\t}_j}{j-n}+\ell^1(n),
\end{aligned}
$$
which yields \er{aa1}, since $\sum_{j\ne n}\frac{\ve_j}{j-n}=\ell^2(n)$
as $n\to \pm \iy$ and ${\rm v.p.}\sum_{j\ne n}\frac1{j-n}=0$.
The proof of other results is similar.\BBox



Below we need results about the Lyapunov function $\D(\l,q)={1\/2}\Tr y(1,\l,q)$ from \cite{K05}.

\begin{lemma}  \lb{T22}
i) The functions $\D (\cdot,\cdot )$ is entire on $\C\ts \mH_\C$ and
the following estimate
\[
|\D (\l,v)|\leq e^{ |\Im \l|+\|v\|_\C}
\]
holds true. Moreover, for each $\l\in\C$ the function $\D(\l,v)$ is even with
respect to $v\in \mH_\C$ and
\[
\D(\l,-v)=\D (\l,v), \ \ \ \ \ v\in \mH_\C.
\]
\end{lemma}

\section{Gauge transformations, even extensions and their properties \lb{Sec3}}
\setcounter{equation}{0}

\subsection{Gauge transformations}

Recall that the operators $\cR, \cF_0, \cF_1,\cF_2$ on $\mH$ are given by
$$
(\cR v)(x)= v(1-x),x\in[0,1],\qq {\rm and}\qq
\cF_0=-\1_2, \ \ \cF_1=J_1\cR,  \qqq      \cF_2=-\cF_1.
$$
In this section, we consider the properties of
ZS-systems under different unitary transformations.
Then the potential $v$ is transformed into another vector,
and corresponding eigenvalues, norming constants also are transformed.
We discuss how the spectral data, eigenfunctions,
and other related quantities move under these transformations.

\begin{lemma}
\lb{Te1} Let $v\in\mH$.
Then the 2-periodic eigenvalues $\l^\pm=(\l^\pm_n)_{n\in\Z}$
and  the fundamental solutions $\vt=\vt(x,\l,v), \vp=\vp(x,\l,v)$ satisfy
\[
\lb{per}
\l^\pm=\l^\pm\circ\cF_j,\ \ \ \ \ j=0,1,2,
\]
\[
\lb{fs0}
\begin{aligned}
\ma\vt_1 &\vp_1 \\ \vt_2&\vp_2\am=\ma\vp_2&-\vt_2\\ -\vp_1&\vt_1\am\circ\cF_0
=\ma\vp_2&\vp_1 \\ \vt_2&\vt_1\am\circ\cF_1=\ma \vt_1& -\vt_2 \\ -\vp_1&\vp_2\am \circ\cF_2\\
\end{aligned}.
\]
Moreover, the Dirichlet eigenvalues $\m=(\mu_n)_{n\in\Z}$,
Neumann eigenvalues $\n=(\nu_n)_{n\in\Z}$ ,
two types mixed eigenvalues $\t=(\t_n)_{n\in\Z}, \vr=(\vr_n)_{n\in\Z}$ and the norming mappings satisfy
\[
\lb{eig}
(\mu,\nu,\t,\vr)=(\nu,\mu,\vr,\t)\circ \cF_0=(\mu,\nu,\vr,\t)\circ\cF_1
=(\nu,\mu,\t,\vr)\circ\cF_2.
\]
\[
\lb{fs1}
\begin{aligned}
(\gr ,\gs ,\gt ,\gu )=(\gs,\gr,\gu,\gt)\circ \cF_0
=(-\gr,-\gs,-\gu,-\gt)\circ \cF_1=(-\gs,-\gr,-\gt,-\gu)\circ \cF_2.
\end{aligned}
\]
\end{lemma}

\no {\bf Proof.} We show the identities \er{per}-\er{fs1} for $\cF_0=-\1_2$.
Let $p=-v$. Lemma  \ref{T22} gives that  $\D(\cdot,-v)=\D(\cdot,v),
$ for all $v\in\mH$
which yields $\l^\pm_n(p)=\l^\pm_n(v)$ for all $n\in\Z$ and the identity \er{per} for $\cF_0$.
The gauge transformation for the system  $Jy'+Vy =\l y$, where $ y=y(x,\l,v)$  gives
$$
\lb{new}
\begin{aligned}
\tes J\big(J{d\/dx}+V\big)J^*f=\tes J f'-Vf=\l f,\qqq {\rm where}
\qq f=Jy(x,\l,v)J^*,
\end{aligned}
$$
since $JV=-VJ$. By the well known results of theory of ODE, we obtain the identity $f=y(x,\l,p)$.
Then it shows that the fundamental solutions $\vt, \vp$ satisfy
$$
y(1,\l,p)=\ma \vt_1& \vp_1\\ \vt_2&\vp_2\am(1,\l,p)=Jy(x,\l,v)J^*=\ma \vp_2& -\vt_2\\ -\vp_1&\vt_1\am(1,\l,v).
$$
This gives the first identity \er{fs0} and
$(\mu,\n,\t,\vr)\circ \cF_0=(\n,\m,\vr,\t),$
since $\m_n,\n_n,\t_n$ and $\vr_n$ are the roots of $\vp_1(1,\l)$, $\vt_2(1,\l)$, $\vp_2(1,\l)$ and $\vt_1(1,\l)$ respectively. Thus we have the identity \er{eig} for $\cF_0$. Moreover, the identity \er{fs0} for $\cF_0$ implies
$$
\begin{aligned}
&e^{-\gr_n(p)}=|\vp_2(1,\m_n(p),p)|=|\vt_1(1,\n_n(v),v)|
=e^{-\gs_n(v)}, \\
&e^{-\gs_n(p)}=|\vt_1(1,\n_n(p),p)|=|\vp_2(1,\m_n(v),v)|
=e^{-\gr_n(v)}, \\
&e^{-\gt_n(p)}=|\vp_1(1,\t_n(p),p)|=|\vt_2(1,\vr_n(v),v)|
=e^{-\gu_n(v)}, \\
&e^{-\gu_n(p)}=|\vt_2(1,\vr_n(p),p)|=|\vp_1(1,\t_n(v),v)|
=e^{-\gt_n(v)}, \\
\end{aligned}
$$
which shows $(\gr,\gs,\gt,\gu)\circ\cF_0=(\gs,\gr,\gu,\gt)$ and the identity \er{fs1} for $\cF_o$ holds true.

We show the identities \er{per}-\er{fs1} for $\cF_1=J_1\cR$. Let $u=\cF_1v=\cR(v_1,-v_2)^\top$.
The gauge transformation for the system  $Jy'+Vy =\l y$, where $ y=y(x,\l,v)$  gives
\[
\lb{new1}
\begin{aligned}
\tes\cF_1\big(J{d\/dx}+V\big)\cF_1^* g=\tes Jg'+V_1g=\l g,
\qqq {\rm where}
\qq g=\cF_1y(x,\l,v)y^{-1}(1,\l,v)\cF_1^*,
\end{aligned}
\]
and $V_1=\cF_1V\cF_1^*=u_1 J_1+u_2 J_2$.
By the known results of ODE theory, we obtain identities
$$
y(x,\l,u)=g(x,\l,v)=J_1y(1-x,\l,v)y^{-1}(1,\l,v)J_1, \qqq y(0,\l,u)=\1_2,
$$
which yields that the fundamental solution $\vt,\vp$ satisfy
$$
y(1,\l,u)=\ma \vt_1& \vp_1\\ \vt_2&\vp_2\am(1,\l,u)=J_1y^{-1}(1,\l,v)J_1=\ma \vp_2& \vp_1\\ \vt_2&\vt_1\am(1,\l,v).
$$
This gives the identity \er{fs0} for $\cF_1$ and
$(\mu,\n,\t,\vr)\circ \cF_1=(\m,\n,\vr,\t),$
since $\m_n,\n_n,\t_n$ and $\vr_n$ are the roots of $\vp_1(1,\l)$, $\vt_2(1,\l)$,
$\vp_2(1,\l)$ and $\vt_1(1,\l)$ respectively. Furthermore,
the identity \er{fs0} and the Wronskian identity
$\vt_1\vp_2-\vt_2\vp_1=1$, imply
\[
\lb{ert}
\left.\begin{aligned}
&e^{-\gr_n(u)}=|\vp_2(1,\m_n(u),u)|
=|\vt_1(1,\m_n(v),v)|=|\vp_2(1,\m_n(v),v)|^{-1}=e^{\gr_n(v)},
\\
&e^{-\gs_n(u)}=|\vt_1(1,\n_n(u),u)|=|\vp_2(1,\n_n(v),v)|
=|\vt_1(1,\n_n(v),v)|^{-1}=e^{\gs_n(v)},
\\
&e^{-\gt_n(u)}=|\vp_1(1,\t_n(u),u)|=|\vp_1(1,\vr_n(v),v)|
=|\vt_2(1,\vr_n(v),v)|^{-1}=e^{\gu_n(v)},
\\
&e^{-\gu_n(u)}=|\vt_2(1,\vr_n(u),u)|=|\vt_2(1,\t_n(v),v)|
=|\vp_1(1,\t_n(v),v)|^{-1}=e^{\gt_n(v)},
\end{aligned}\right.
\]
for all $n\in\Z$. Due to \er{ert} we have $\gr_n\circ\cF_1=-\gr_n, \gt_n\circ\cF_1=-\gu_n$.
Since the operators $\cF_0$ and $\cF_1$ are commute, i.e.
$\cF_0\cF_1=\cF_1\cF_0$, we obtain for all $n\in\Z$,
$$
\gs_n\circ\cF_1=(\gr_n\circ\cF_0)\circ\cF_1=(\gr_n\circ\cF_1)\circ\cF_0
=-\gr_n\circ\cF_0=-\gs_n, \qq and \qq \gu_n\circ\cF_1=-\gt_n.
$$
Thus the identities for $\cF_1$ have been proved. Moreover, we have
$$
\begin{aligned}
&\gr_n\circ(\cF_0\cF_1)=\gs_n\circ\cF_1=-\gs_n, \qqq
\gs_n\circ(\cF_0\cF_1)=\gr_n\circ\cF_1=-\gr_n,  \\
&\gt_n\circ(\cF_0\cF_1)=\gu_n\circ\cF_1=-\gt_n, \qqq
\gu_n\circ(\cF_0\cF_1)=\gt_n\circ\cF_1=-\gu_n.
\end{aligned}$$
The identities \er{per}-\er{fs1} for $\cF_0, \cF_1$ and $\cF_2=\cF_0\cF_1$ yield \er{per}-\er{fs1} for $\cF_2$.       \BBox

The next lemma shows how the normalizing mappings behave under the transformation.

\begin{lemma} Let $\ga_n ^\bu, \gb_n^\bu,\gc_n^\bu, \gd_n^\bu, n\in\Z$ be given by \er{ss1}, \er{ss2}.
The normalizing mappings $\ga=(\ga_n )_{n\in\Z}, \gb=(\gb_n )_{n\in\Z}$,
$\gc=(\gc_n )_{n\in\Z}, \gd=(\gd_n )_{n\in\Z},$ defined by \er{mnc} satisfy
\[
\lb{NC0}
(\ga,\gb,\gc,\gd)=(\gb,\ga,\gd,\gc)\circ \cF_0
=(\mathring\ga,\mathring\gb,\mathring\gd,\mathring\gc)\circ \cF_1 = (\mathring\gb,\mathring\ga,\mathring\gc,\mathring\gd)\circ \cF_2,
\]
where $\mathring\ga=(\mathring\ga_n )_{n\in\Z}, \cdots$ are given by
$$
\mathring\ga_n =\ga_n-2\gr_n , \qq \mathring\gb_n= \gb_n-2\gs_n, \qq
\mathring\gc_n=\gc_n-2\gt_n , \qq \mathring\gd_n=\gd_n-2\gu_n .
$$

\end{lemma}

\no {\bf{Proof.}} We prove the identity for $\cF_0$ in \er{NC0}, the proof for $\cF_1$ and $\cF_2$ is similar. From \er{ss1} and \er{ss2}, we obtain that the components of normalizing mappings $\ga=(\ga_n )_{n\in\Z}, \gb=(\gb_n )_{n\in\Z}$, $\gc=(\gc_n )_{n\in\Z}, \gd=(\gd_n )_{n\in\Z}$ have the following form:
$$
\ga_n=\gr_n -\ga_n ^\bu, \qq \gb_n=\gs_n -\gb_n^\bu, \qq
\gc_n=\gt_n -\gc_n^\bu, \qq  \gd_n=\gu_n -\gd_n^\bu.
$$
Let $p=\cF_o v$. From \er{fs0}, we have
$$
(-\dot\vp_1,\dot\vp_2, \dot\vt_1,-\dot\vt_2)(1,\cdot,p)=(\dot\vt_2,\dot\vt_1,
\dot\vp_2,\dot\vp_1)(1,\cdot,v).
$$
Combining the previous identity with \er{eig}, we deduce that for all $n\in\Z$
\[
\lb{bb}
\begin{aligned}
&\ga_n^\bu(p)=\ln|\dot\vp_1(1,\m_n(p),p)|
=\ln|\dot\vt_2(1,\n_n(v),v)|=\gb_n^\bu(v), \\
&\gb_n^\bu(p)=\ln|\dot\vt_2(1,\n_n(p),p)|
=\ln|\dot\vp_1(1,\m_n(v),v)|=\ga_n^\bu(v), \\
&\gc_n^\bu(p)=\ln|\dot\vp_2(1,\t_n(p),p)|
=\ln|\dot\vt_1(1,\vr_n(v),v)|=\gd_n^\bu(v), \\
&\gd_n^\bu(p)=\ln|\dot\vt_1(1,\vr_n(p),p)|
=\ln|\dot\vp_2(1,\t_n(v),v)|=\gc_n^\bu(v).
\end{aligned}
\]
Applying $\cF_0$ to the normalizing constants $\ga_n , \gb_n, \gc_n, \gd_n, n\in\Z$ we have
\[
\lb{NCm}
\begin{aligned}
\ga_n \circ \cF_0&=\gs_n -\gb_n^\bu=\gb_n, \qqq \gb_n\circ \cF_0=\gr_n -\ga_n ^\bu=\ga_n ,\\
\gc_n\circ \cF_0&=\gu_n -\gd_n^\bu=\gd_n, \qqq \gd_n\circ \cF_0=\gt_n -\gc_n^\bu=\gc_n,
\end{aligned}\]
since \er{fs1} for $\cF_o$ and \er{bb} hold true. The first identity of \er{NC0} follows.

\no The application of $\cF_1$ to the normalizing constants $\ga_n ^\bu,
\gb_n^\bu,\gc_n^\bu$, $\gd_n^\bu, n\in\Z$, where $u=\cF_1v$, yields
$$
\begin{aligned}
&\ga_n^\bu(u)=\ln|\dot\vp_1(1,\m_n(u),u)|
=\ln|\dot\vp_1(1,\m_n(v),v)|=\ga_n^\bu(v), \\
&\gb_n^\bu(u)=\ln|\dot\vt_2(1,\n_n(u),u)|
=\ln|\dot\vt_2(1,\n_n(v),v)|=\gb_n^\bu(v),
\\
&\gc_n^\bu(u)=\ln|\dot\vp_2(1,\t_n(u),u)|
=\ln|\dot\vt_1(1,\vr_n(v),v)|=\gd_n^\bu(v),
\\
&\gd_n^\bu(u)=\ln|\dot\vt_1(1,\vr_n(u),u)|
=\ln|\dot\vp_2(1,\t_n(v),v)|=\gc_n^\bu(v).
\end{aligned}
$$
These, together with \er{fs1} for $\cF_1$, give
\[
\lb{NCm1}
\begin{aligned}
\ga_n \circ \cF_1&=-\gr_n -\ga_n^\bu=\ga_n-2\gr_n, \qqq
\gb_n\circ \cF_1=-\gs_n -\gb_n ^\bu=\gb_n-2\gs_n, \\
\gc_n\circ \cF_1&=-\gu_n -\gd_n^\bu=\gd_n-2\gu_n, \qqq
\gd_n\circ \cF_1=-\gt_n -\gc_n^\bu=\gc_n-2\gt_n.
\end{aligned}
\]
Then the second identity of \er{NC0} holds true. Combining \er{bb}-\er{NCm1} and
$\cF_2=\cF_0\cF_1$, we obtain
$$
\begin{aligned}
\ga_n\circ \cF_2&=-\gs_n -\gb_n ^\bu=\gb_n-2\gs_n, \qqq
\gb_n\circ \cF_2=-\gr_n -\ga_n^\bu=\ga_n-2\gr_n,\\
\gc_n\circ \cF_2&=-\gt_n -\gc_n^\bu=\gc_n-2\gt_n, \qqq
\gd_n\circ \cF_2=-\gu_n -\gd_n^\bu=\gd_n-2\gu_n.
\end{aligned}
$$
Thus the  identity \er{NC0} has been proved.     \BBox

\subsection{Shifting mappings}
Applying a ratation on the plane $\cF =e^{{\pi\/2}x J}$, we obtain
\[
\lb{eq2}
\cF \big(J\tes{d\/dx}+V\big)\cF ^*=J\tes{d\/dx}+{\pi\/2}+V_u,
\qq V_u:=\ma u_1 & u_2\\ u_2 &-u_1\am, \qq \gF =\tes e^{{\pi}x J},
\]
where
$$
V_u=\cF  V\cF ^*=\gF V=\ma v_1c+v_2s & -v_1s+v_2c \\ v_2c-v_1s & -(v_1c+v_2s)\am,
$$
and the vector $u$ is given by
$$
u=\ma u_1 \\ u_2\am=e^{\pi xJ}v=\ma v_1c+v_2s \\ -v_1s+v_2c \am, \qqq
\ca c=\cos\pi x \\ s=\sin\pi x \ac .
$$
We consider relations between the spectra of the operator before and after its transformation.

\begin{lemma}
\lb{Te3}
Let $\gF =\tes e^{{\pi}x J}$ and the shifted operator $\cS:\gJ^o\to \gJ^1$ be given by
$
 (\cS z)_n=\tes z_n-{\pi\/2}$ for all $z=(z_n)_{n\in\Z}\in\gJ^o$.
Then

\no i) The shifted operator $\cS:\gJ^o\to\gJ^1 $ is a bijection between $\gJ^o$ and $\gJ^1$.

\no ii) The solutions $\vt,\vp$  satisfy
\[
\lb{Ca1}
\tes
\ma \vt_1 & \vp_1 \\ \vt_2 & \vp_2\am (x,\l,v)
=\ma -\vt_2 & -\vp_2 \\ \vt_1 & \vp_1\am (x,\l-{\pi\/2}, \gF (v)).
\]
Moreover, the 1-spectra mappings $\m, \n, \t,\vr$, norming mappings
$\gr, \gs, \gt,\gu$ and the normalizing mappings $\ga, \gb, \gc, \gd$  satisfy
\[
\lb{cA}
(\t,\vr)\circ \gF=(\cS\m,\cS\n), \qqq (\gt,\gu)\circ \gF=(\gr,\gs), \qqq
(\gc,\gd)\circ \gF=(\ga,\gb).
\]

\end{lemma}

\no {\bf{Proof.}} i) We show that the operator $\cS$ is an injection.
In $\cS \a=\cS \b$ for  sequences
$\a=(\a_n)_{n\in\Z}, \b=(\b_n)_{n\in\Z}\in\gJ^o$, then their components satisfy $\tes \a_n=\b_n,$ for all $n\in\Z$.

We show a surjection.
Let $\b=(\b_n)_{n\in\Z}\in\gJ^1$. Define the sequence $\a=(\a_n)_{n\in \Z}$ by $\a_n=\tes\b_n+ \frac{\pi}{2}$. Then $\a\in\gJ^o$, and a direct computation shows that $\cS\a =\b$. Therefore, the operator $\cS$ is a bijection between $\gJ^o$ and $\gJ^1$.

\no ii) Let $\cF=e^{{\pi\/2} x J}$ and $V_u$ be defined by \er{eq2}.
Due to \er{eq2} we deduce that solution $f=y(x,\z,u)$ of the system $Jf'+V_uf=\z f$,
where $u=\gF v$ has the form
$$
y(x,\z,u)=\cF (x) y(x,\l,v), \qqq \tes\z=\l-{\pi\/2}, \qqq y(0,\z,u)=\1_2.
$$
Hence, due to $\cF(1)=J$ the fundamental solutions $\vt, \vp$ satisfy
\[
\lb{fund}
y(1,\z,u)=\ma \vt_1 & \vp_1 \\ \vt_2 & \vp_2\am(1,\z,u)
=J y(1,\l,v)=\ma \vt_2 & \vp_2 \\ -\vt_1 & -\vp_1\am(1,\l,v),
\]
which implies \er{Ca1}. Recall that
$\m_n, \n_n, \t_n,\vr_n$ are the roots of
$\vp_1(1,\l), \vt_2(1,\l)$, $\vp_2(1,\l)$ and $\vt_1(1,\l)$ respectively.
Identity \er{fund} implies $
\tes\t_n(u)=\m_n(v)-{\pi\/2}, \ \vr_n(u)=\n_n(v)-{\pi\/2},
$
and
$$
\begin{aligned}
\tes
e^{-\gt_n(u)}=|\vp_1(1,\t_n(u),u)|=|\vp_2(1,\t_n(u)+{\pi\/2},v)|
=|\vp_2(1,\m_n(v),v)| =e^{-\gr_n(v)},
\\
\tes e^{-\gu_n(u)}=|\vt_2(1,\vr_n(u),u)|=|\vt_1(1,\vr_n(u)+{\pi\/2},v)|
=|\vt_1(1,\n_n(v),v)| =e^{-\gs_n(v)},
\end{aligned}
$$
 for all $n\in\Z$, which yields the first and second identities of \er{cA}.

From \er{ss1} and \er{ss2}, we obtain that the components of normalizing
mappings $\ga=(\ga_n )_{n\in\Z}, \gb=(\gb_n )_{n\in\Z}$, $\gc=(\gc_n )_{n\in\Z},
\gd=(\gd_n )_{n\in\Z}$ have the following form:
$$
\ga_n=\gr_n -\ga_n ^\bu, \qq \gb_n=\gs_n -\gb_n^\bu, \qq
\gc_n=\gt_n -\gc_n^\bu, \qq  \gd_n=\gu_n -\gd_n^\bu.
$$
From \er{fund}, we obtain
$
\begin{aligned}
\tes(\dot\vt_1,\dot\vp_2)(1,\l-{\pi\/2},u)=(\dot\vt_2,-\dot\vp_1)
(1,\l,v),
\end{aligned}$
which gives
$$
\begin{aligned}
\gc_n^\bu(u)=\ln|\dot\vp_2(1,\t_n(u),u)|=\ln|\dot\vp_1(1,\m_n(v),v)|
=\ga_n^\bu, \\
\gd_n^\bu(u)=\ln|\dot\vt_1(1,\vr_n(u),u)|=\ln|\dot\vt_2(1,\n_n(v),v)|
=\gb_n^\bu.
\end{aligned}$$
This, jointly  with the second identity of \er{cA}, implies
$\gc_n\circ \gF=\gr_n -\ga_n^\bu=\ga_n $ and
$\gd_n\circ \gF=\gs_n -\gb_n^\bu=\gb_n,$
which yields the third identity of \er{cA}.      \BBox

\subsection{Even extensions}

We define the space $\wt\mH=L^2([0,2],\R)^2$ and the even-odd spaces
by
$$
\wt\mH_{eo}=\big\{v\in \wt\mH: v_1(2-x)=v_1(x), \qq v_2(2-x)=-v_2(x), \qq  x\in(0,2)\big\}.
$$
For the vector $v\in\mH$ we define
the even-odd extension $\cE :v\to \wt v$ acting from $\mH$ into $\wt\mH_{eo}$ and the corresponding matrix  $\wt V$ by
\[
\lb{TV}
 v=\ma  v_1\\ v_2\am\to
\cE v=\wt v=\ma \wt v_1\\ \wt v_2\am,\qqq \qq
\wt V=\ma \wt v_1&\wt v_2\\ \wt v_2&-\wt v_1\am,
\]
where
\[
\lb{q}
\wt v(x)=v(x),\qq 0<x<1,\qqq {\rm and}\qq
\wt v(x)=J_1v(2-x),\qq 1<x<2.
%
%
%
\]
For the matrix $\wt V$ given by \er{TV}
we introduce an operator $\wt T_{per}y=J y'+\wt Vy$ with 4-periodic boundary conditions. Let $\wt \l_n^\pm, n\in \Z$ be eigenvalues of $\wt T_{per}$ labeled by the standard way $\cdots<\wt \l^-_{n-1}\le\wt \l^+_{n-1}<\wt \l^-_{n}\le\wt \l^+_{n}<\cdots$,
where the equality $\wt\l_n^-=\wt\l_n^+$ means that $\wt\l_n^-$ is a double eigenvalue.
Here $\wt \l^\pm_{2n}$ is an eigenvalue with  2-periodic boundary conditions
and $\wt \l^\pm_{2n+1}$ is an eigenvalue with  anti 2-periodic boundary conditions. The eigenvalues $\wt \l_n^\pm$ have asymptotics
$$
\wt \l^\pm_n=\tes{n\pi\/2}+o(1) \qqq \text{as} \ \ \ n\to\pm\iy.
$$
Define a gaps $\wt \g_n=(\wt \l_n^-,\wt \l_n^+)$
with the length $|\wt \g_n|\ge0$.
We consider the ZS-systems on the interval $[0,2]$ under Dirichlet and Neumann
  boundary conditions:
\begin{equation}
\label{DNe} Jf'+\wt Vf=\l f,\qqq
\begin{aligned}
&     \qqq f_1(0)= \textstyle f_1(2)=0, \qq \{\m_n, n\in \Z\} \ Dirichlet
\\
&   \qqq  f_2(0)=\textstyle f_2(2)=0,  \qq \{\n_n, n\in \Z\} \ Neumann
\end{aligned}\ \ \ .
\end{equation}
Let $\wt \mu_n=\m_n(\wt v)$ and $\wt \nu_n=\n_n(\wt v), n\in\Z$ be the Dirichlet and Neumann eigenvalues respectively.
The next lemma shows their positions.
Recall that
$$
\mH_{eo}=\big\{v\in \mH: v=J_1\cR v\big\},\qq \mH_{oe}=\big\{v\in \mH: v=-J_1\cR v\big\}.
$$

\begin{lemma}
\lb{Te4}i) The mapping $\gF= e^{{\pi}x J}:\mH_{eo}\to\mH_{oe}$ is a bijection between $\mH_{eo}$ and $\mH_{oe}$.

\no ii) Let $v\in\mH_{eo}$. Then the Lyapunov function $\D$
and the norming mappings satisfy
\[
\lb{Te1-1}
\D=\vp_2(1,\cdot)=\vt_1(1,\cdot),
\]
\[
\lb{ee1}
\gr(v)=\gs(v)=\gt(\gF v)=\gu(\gF v)=0\qq \forall \  v\in\mH_{eo}.
 \]
If $\l\in \R$, then
\[
\lb{Te1-2}
\D^2(\l)=1  \qqq \Leftrightarrow    \qqq  \vp_1(1,\l)=0 \qq  \ or \qq \vt_2(1,\l)=0.
\]
Moreover, the Dirichlet and Neumann eigenvalues form the endpoints of gaps $\g_n$:
\[
\lb{GN}
\g_n=(\mu_n,\nu_n) \qqq \text{or} \qqq \g_n=(\nu_n,\mu_n), \qq \forall \  n\in\Z.
\]
\end{lemma}

\no {\bf Proof.} i) Let $\gF=e^{{\pi}x J}$.
For any $v\in\mH_{eo}$ we have
$$
u=\ma u_1 \\ u_2 \am=\gF v=\ma v_1c+v_2s \\ -v_1s+v_2c \am, \qqq
\ca c=c(x)=\cos\pi x \\ s=s(x)=\sin\pi x \ac,
$$
where  $u_1=v_1c+v_2s$ is  odd and $u_2=-v_1s+v_2c$ is even, since $(c,s)^\top\in \mH_{oe}$, which yields $u\in\mH_{oe}$.
Similar arguments give that for any $u\in\mH_{oe}$ there exists a unique
$v\in\mH_{eo}$ such that $\gF v=u$.

We show \er{ee1}. Let $v\in\mH_{eo}$. From \er{Te1-1} and the results of $\D$, we obtain
$
\D(\m_n)=\vp_2(1,\m_n,v)=(-1)^n$ for all $n\in\Z,  v\in\mH_{eo}$,
which yields $\gr_n(v)=0$. This, together with \er{cA},
gives $\gt_n\circ\gF=\gr_n=0$, $n\in\Z$. Similarly  we have
$
\D(\n_n)=\vt_1(1,\n_n,v)=(-1)^n$ for all $n\in\Z$.
This identity and \er{cA} give $\gu_n\circ\gF=\gs_n=0$, $n\in\Z$.

\no ii) Recall that $\cF_1=J_1\cR$. Since $v\in\mH_{eo}$,
the identity $\cF_1 v=v$ holds true. Then the second identity of \er{fs0} shows
$$
\ma \vp_2 & \vp_1 \\ \vt_2 & \vt_1 \am(x,\l,\cF_1v)=\ma \vp_2 & \vp_1 \\ \vt_2 & \vt_1 \am(x,\l,v)
=\ma \vt_1 & \vp_1 \\ \vt_2 & \vp_2 \am(x,\l,v), \qq (x,\l)\in[0,1]\times\R,
$$
which yields \er{Te1-1}. Substituting \er{Te1-1} into the
Wronskian $\vt_1\vp_2-\vt_2\vp_1=1$ we get $\D^2(\l)-1=\vt_2(1,\l)\vp_1(1,\l)$,
this gives \er{Te1-2} and \er{GN}.
 \BBox

\begin{lemma} \lb{Tet}
Let $\cE:\mH\to\wt\mH_{eo}$ be an even-odd extension.
Then the Dirichlet  $\wt \m_{n}$ and the Neumann eigenvalues $\wt \n_n,$ and 4-periodic eigenvalues $\wt \l_n^\pm$ satisfy
\[
\lb{Ev1}
(\wt \m_{2n-1},\wt \m_{2n},\wt \n_{2n-1},\wt \n_{2n})\circ\cE=(\t_n,\m_n,\vr_n,\n_n),
\]
\[
\lb{Ev4}
\left\{
\begin{aligned}
&\{\wt \l_{2n-1}^-,\wt \l_{2n-1}^+\}\circ\cE=\{\vr_{n},\t_{n}\}\circ\cF_1
=\{\vr_{n},\t_{n}\}\circ\cF_2 \\
&\{\wt \l_{2n}^-,\wt \l_{2n}^+\}\circ\cE=\{\m_{n},\n_{n}\}\circ\cF_1
=\{\m_{n},\n_{n}\}\circ\cF_2
\end{aligned}
\right.   ,
\]
for all $n\in\Z$, where $\{\x,\z\}$ denotes a set of two points $\x, \z\in \R$.
\end{lemma}

\no {\bf{Proof.}} Let $v\in\mH$. Let $\m_n$ and $\p_n,~n\in\Z$ be the Dirichlet eigenvalues and the corresponding eigenfunctions. We have $J\p_n'+V\p_n=\mu_n\p_n$
and $\p_{n1}(0)=\p_{n1}(1)=0.$ Define a function
$$
\tes
f(x)=
\left\{\begin{aligned}&\p_n(x), \qqq &&x\in(0,1)
\\ &-J_1\p_n(2-x),&& x\in(1,2)
\end{aligned}\right. .
$$
It is continuous on at $x=1$, since
$$
\tes
f(1+0)=-J_1\p_n(1)=\ma 0 \\ \p_{n2}(1) \am=\p_n(1)=f(1-0).
$$
Let $\wt v=\cE v\in\wt\mH_{eo}$, and the corresponding matrix $\wt V$
be given by \er{TV}. The direct calculation
shows that $Jf'+\wt Vf=\mu_n f$
on $[0,2]$ with $f_1(0)=f_1(2)=0$, which implies that $\m_n$
are the Dirichlet eigenvalues for $\wt v\in\wt\mH_{eo}$.
Lemma \ref{Te4} gives that Dirichlet eigenvalues are located at the end of the gaps.
Then by the basic asymptotics of eigenvalues,
we obtain that $\mu_n$ coincide with the Dirichlet eigenvalue $\wt\m_{2n}$ for $|n|\to\iy$.

Let $\t_n$ and $\phi_n$, $n\in\Z$, be the mixed eigenvalues and the corresponding eigenfunctions,
which satisfies $\phi_{n1}(0) =\phi_{n2}(1)=0$. Define a function
$$
g(x)=
\left\{\begin{aligned}&\phi_n(x), \qqq &&x\in(0,1)
\\ &J_1\phi_n(2-x),&& x\in(1,2)
\end{aligned}\right.   .
$$
It is continuous on at $x=1$, since
$$
g(1+0)=J_1\phi_n(1)=\ma \phi_{n1}(1) \\ 0 \am=\phi_n(1)=g(1-0).
$$
The direct calculation gives that $g$ satisfies $Jg'+\wt Vg=\t_n g$ on $[0,2]$
and $g_1(0)=g_1(2)=0$, which implies that $\t_n$
are the Dirichlet eigenvalues for $\wt v\in\wt\mH_{eo}$.
Lemma \ref{Te4} gives that Dirichlet eigenvalues are located at the end of the gaps.
Then by the basic asymptotics of eigenvalues,
we obtain that $\t_n$ coincide with the Dirichlet eigenvalue $\wt\m_{2n-1}$ for $|n|\to\iy$.

Let $\n_n$ be Neumann eigenvalues  and $\x_n$ be the corresponding eigenfunctions such that $
J\x_n'+V\x_n=\n_n\x_n$ and $\x_{n2}(0)=\x_{n2}(1)=0$ for all $n\in\Z.
$
Define a function
$$
h(x)=
\left\{\begin{aligned}&\x_n(x), \qqq &&x\in(0,1)
\\ &J_1\x_n(2-x),&& x\in(1,2)
\end{aligned}\right.   .
$$
It is continuous on at $x=1$, since
$$
h(1+0)=J_1\x_n(1)=\ma \x_{n1}(1) \\ 0 \am=\x_n(1)=h(1-0).
$$
The direct calculation gives that $h$ satisfies $Jh'+\wt Vh=\t_n h$ on $[0,2]$
and $h_2(0)=h_2(2)=0$, which implies that $\n_n$
are the Neumann eigenvalues for $\wt v\in\wt\mH_{eo}$.
Lemma \ref{Te4} gives that Neumann eigenvalues are located at the end of the gaps.
Then by the basic asymptotics of eigenvalues,
we obtain that $\n_n$ coincide with the Neumann eigenvalue $\wt\n_{2n}$ for $|n|\to\iy$.

Let the mixed eigenvalues $\vr_n$ and the corresponding  eigenfunctions $\e_n$ satisfy
$
J\e_n'+V\e_n=\vr_n\e_n$ and $\e_{n2}(0)=\e_{n1}(1)=0$ $ n\in\Z.
$
Define a function
$$
G(x)=
\left\{\begin{aligned}&\e_n(x), \qqq &&x\in(0,1)
\\ &-J_1\e_n(2-x),&& x\in(1,2)
\end{aligned}\right.  .
$$
It satisfies
$$
G(1+0)=-J_1\e_n(1)=\ma 0 \\ \e_{n2}(1) \am=\e_n(1)=G(1-0).
$$
The direct calculation gives that $G$ satisfies $J{d\/dx}G+\wt VG=\t_n G$ on $[0,2]$
and $G_2(0)=G_2(2)=0$, which implies that $\vr_n$
are the Neumann eigenvalues for $\wt v\in\wt\mH_{eo}$.
Lemma \ref{Te4} gives that Neumann eigenvalues are located at the end of the gaps.
Then by the basic asymptotics of eigenvalues,
we obtain that $\vr_n$ coincide with the Neumann eigenvalue $\wt\n_{2n-1}$ for $|n|\to\iy$.
Thus we obtain \er{Ev1}. Applying Lemma \ref{Te1} to \er{GN} we obtain \er{Ev4}.       \BBox

\begin{corollary}
\label{T2smb}

\no i) All 2-spectra mappings $\t\star \m, \varrho\star \n$,  $\varrho\star
\m$ and $\t\star \n$   are isomorphic
and satisfy
\begin{equation}
\label{enm4}
\begin{aligned}
\t\star \m=(\varrho\star \n)\circ \cF_o=(\varrho\star \m)\circ \cF_1=(\t\star
\n)\circ \cF_2.
\end{aligned}
\end{equation}
\no ii) All mappings $\m\ts \gr, ~\n \ts \gs $,
$\m \ts (-\gr)$, $\n \ts (-\gs) $, $\t\ts \gt,~ \vr \ts \gu $, $\vr \ts (-\gu) $
and $\t \ts (-\gt)$ are isomorphic and satisfy
\[
\label{enm14}
\begin{aligned}
\m\ts \gr =(\n \ts \gs )\circ \cF_o
=(\m \ts (-\gr) )\circ \cF_1=(\n \ts (-\gs))\circ \cF_2,
\end{aligned}
\]
\[
\label{enm14x}
\t\ts \gt =(\vr \ts \gu )\circ \cF_o
=(\vr \ts (-\gu) )\circ \cF_1=(\t \ts (-\gt))\circ \cF_2.
\]

\no iii) Let  $\cE  : \mH\to\wt\mH_{eo}$ be  the extension  given by \er{q}. Let $\wt \m_{n}, \wt \n_{n}, n\in\Z$ be the Dirichlet and Neumann eigenvalues  for $\wt v=\cE v\in\wt\mH_{eo}$ respectively. Then the mappings $v\to \wt \m=(\wt \m_{n})_{n\in \Z}$ and $v\to \wt \n=(\wt \n_{n})_{n\in \Z}$  satisfy
\[
\label{dwtd}
\t\star \m=\wt \m\circ \cE,\qqq \vr\star \n=\wt \n\circ \cE.
\]

\no iv) Let $\gF$ and $\cS$ be defined by \er{eq2}. The mappings $(\cS\m)\ts \gr$, $(\cS\n)\ts \gs$
are isomorphic and satisfy
\begin{equation}
\label{cA4}  (\cS\m)\ts \gr=(\t \ts \gt)\circ \gF,\qqq (\cS\n)\ts \gs=(\vr \ts \gu)\circ \gF.
\end{equation}
\end{corollary}

\no {\bf Proof.} i) From Lemma \ref{Te1} we deduce that all 2-spectra mappings $\t\star \m, \varrho\star \n$,  $\varrho\star
\m$ and $\t\star \n$   are isomorphic and satisfy \er{enm4}.

\no ii)  From Lemma \ref{Te1} we deduce that all 2-spectra mappings
$\m\ts \gr, ~\n \ts \gs $, $\m \ts (-\gr)$, $\n \ts (-\gs) $,
$\t\ts \gt,~ \vr \ts \gu $, $\vr \ts (-\gu) $
and $\t \ts (-\gt)$ are isomorphic and satisfy \er{enm14}.

\no iii) From \er{Ev1} we obtain the identity $\t\star \m=\wt \m\circ \cE$ in \er{dwtd}.  The proof for $\vr\star \n$
is similar.
\BBox


\section {Isomorphic inverse problems on the finite interval \lb{Sec4}}
\setcounter{equation}{0}

\subsection {Preliminary results}

We recall the well known results about analytic functions in the Hilbert space, see p. 138 \cite{PT87}.

\begin{theorem}\lb{TAi2}
Let $f: \mD \to \cH$ be a map from an open subset $\mD$ of a complex
Hilbert space $\mH$ into a Hilbert space $\cH$ with orthonormal basis $e_n,
n\in\Z$. Then $f$ is analytic on $\mD$ if and only if is locally
bounded, and each "coordinate function"
$
f_n=\langle f,e_n \rangle :\mD\to C
$
is analytic on $\mD$. Moreover, the derivative of $f$ is given by the
derivatives of its "coordinate functions":
$$
f'(v)h=\sum_{n\in\Z}\langle f'(v)h,e_n \rangle,\qq h\in\mH.
$$

\end{theorem}

Below we need following results about  basises from p. 163  \cite{PT87}.

\begin{theorem}
\label{Tfa1} Let $e_n^o, n\in\Z$ be an orthogonal basis of the Hilbert space
$\mH$. Suppose  $e_n, n\in\Z$ is another sequence of vectors in $\mH$ that either spans or is linear independent. If, in addition,
$$
\sum_{n\in\Z}\|e_n-e_n^o\|^2<\iy,
$$
then $e_n, n\in\Z$ is also a basis of $\mH$. Moreover, the map
$v\to (\lan v,e_n\ran )_{n\in \Z}$  is a linear isomorphism between $\mH$ and $\ell^2$.
\end{theorem}


Recall results  about the transformations
both for frozen  norming constants and frozen eigenvalues  from \cite[Th~3.1,3.2]{DK00} about "locally free perameters".
In the case of Schr\"odiger operators results about "locally free perameters"
were decribe in the book [PT87], see p. 91, 111.

 \begin{theorem}
\label{Tfp}
Let $\n(v)\ts \gs(v)\in \cJ^o\ts \ell^2$ for some $v\in \mH $ and let $m\in \Z$. Then

\no i) For any 
 sequence $\x=(\x_n)_{n\in\Z}\in \cJ^o$, where
$\x_n=\n_n(v), n\ne m$ and $\x_m\in (\n_{m-1}(v),\n_{m+1}(v) )$
there exists a potential $w\in \mH$ such that $(\n\ts \gs)
(w)=\x\ts \gs(v)$.
Moreover, if $v\in \mH_{eo}$, then  $\gs(v)=0$ and $w\in \mH_{eo}$.

\no ii) For any  sequence   $\gl=(\gl_n)_{n\in\Z}\in \ell^2$, where
$\gl=\gs_n(v), n\ne m$ and $\gl_m\in\R$ there exists  $w\in \mH$ such that $(\n\ts \gs)(w)=(\n(v), \gl)$.

\end{theorem}

This theorem shows that spectral data $(\n_n, \gs_n)_{n\in \Z}$
are "locally free parameters". It means the following:
we fix all parametrs  except one. The last parameter can be moved to any point in the interval. This interval is  finite in
the case of the eigenvalue  and this interval is the real line  in
the case of the norming constant.  For each new parameter there exists a potential from $\mH$ such that .

We reformulate results of P\"oschel and Trubowitz \cite{PT87} for the ZS-systems.

 \begin{theorem}
\label{Taa} The mapping $f=\m\ts \gr: \mH\to \gJ^o\ts \ell^2$
 has the following properties:

\no i) The mapping $f$ is real analytic.

\no ii) 
The mapping $f$
is a real analytic local isomorphism between $\mH$ and $\cJ^o\ts \ell^2$.

\no iii) The mapping $f$ is one-to-one.

\no iv)  The mapping $f$ is a surjection.

Moreover,  $f$ is a RAB between $\mH$ and $\gJ^o\ts \ell^2$.
\end{theorem}

\no {\bf Proof.} Properties i) and ii) are proved in Theorem \ref{Tdnc}.

\no iii) The injection of the mapping for the ZS-systems is well known fact, see e.g., \cite{A14}, \cite{GL66}.

\no
iv)  We show a surjection. Let $\f=(\f_n)_{n\in\Z}\in \cJ^o\ts \ell^2$  and
$\f_n=(\e_n,\z_n) $.  Consider the cut sequence $\f^m=(\f_n^m)_{n\in \Z}$, where
$$
\begin{aligned}
\f_n^m=f_n(0),\ \ \forall \  |n|\le m, \qqq \f_n^m=(\f_n) \qq \forall \  |n|>m.
\end{aligned}
$$
The sequence  $\f^m$ converges to $f(0)=(\m(0),\gr(0))$ as $m\to\iy$. Thus they must be contained in the open image of the map $f$. Then we have $f(u^m)=\f^m$ for $m$ large enough.
It remains to shift the first $2m+1$ eigenvalues  $\m_n(u^m)=\m_n^o$ to $\e_n$ and both $2m+1$ norming constants of $\gr_n(u^m)=0$ to $\z_n, |n|\le m$. We can do it via Theorem \ref{Tfp}, changing only finite  number of spectral parameters. Here we use the proof from the great book \cite{PT87},
 see also \cite{DK00}.

Properties i)-iv) imply that $f$ is a RAB between $\mH$ and $\gJ^o\ts \ell^2$.
 \BBox.

 \begin{theorem}
\label{Tde}
 The 1-spectra mapping $\m: \mH_{eo}\to \gJ^o$
 has the following properties:

\no i) The mapping $\m$ is real analytic.

\no ii) The mapping $\m$ is a real analytic local isomorphism.

\no iii) The mapping $\m$ is one-to-one.

\no iv)  The mapping $\m$ is a surjection.

Moreover,  $\m$ is a RAB between $\mH_{eo}$ and $\gJ^o$.

\end{theorem}
\no {\bf Proof.}
i) Lemma \ref{T23} gives that the mapping $v\to \m$ acting from $\mH_{eo}$ into $\gJ^o$ is real analytic.

ii) We show that the mapping $v\to \m$ is the local real analytic isomophism.
Let $\m_n'(v)={\pa \m_n(v)\/\pa v}$ for shortness.
By Lemma \ref{T25} and Theorem \ref{T25}, for fixed $v\in\mH_{eo}$
the operator $\m'(0)$ is the Fourier transformation from $\mH_{eo}$ onto $\ell^2$ and the operator $ \m'(v)- \m'(0)$ is a compact. Thus $\m'(v)$ is a Fredholm operator.

We prove that the operator $\m'(v)$ is invertible. Assume that it is not
invertible. Then there exists $h\in \mH_{eo}, h\ne 0$, which is a solution of the equation
$$
\m'(v)h=0\ \ \ \ \Leftrightarrow \ \ \ \ \ \   \ \lan \m_n'(v),h\ran=0, \ \ \forall \ n\in \Z.
$$
Due to Lemma \ref{T25},
the sequence $\m_n'(v), n\in \Z$ is linearly independent.
Then  using (2.38) and Theorem \ref{Tfa1} we deduce that the sequence $( \m_n')_{ n\in \Z}$  forms a basis of $\mH_{eo}$.
This implies that $h=0$, since the sequence $(\m_n'(v))_{ n\in \Z}$  forms a basis of $\mH_{eo}$. Thus due to the Invese Function Theorem, the operator $\m'(v)$ is invertible and $v\to \m$ is a real analytic local isomorphism.

\no iii) and iv) The proof repeats the case of Theorem \ref{Taa}.

 Thus  the mapping $\m:\mH_{eo}\to \gJ^o$ is a RAB   between $\mH_{eo}$ and $\gJ^o$.
\BBox

We reformulate Theorem \ref{Taa} for the shifting mapping.

 \begin{theorem}
\label{Tsa}
i)  The mappings $(S\m)\ts \gr $ and $\t \ts \gt$ acting from $\mH$ to $\gJ^1\ts \ell^2$ are  a RAB between $\mH$ to $\gJ^1\ts \ell^2$ and
satisfy $(\cS\m)\ts \gr=(\t \ts \gt)\circ \gF$.

\no ii) Each of the mappings $S\m$ and $\t$ acting from $\mH_{eo}$ into $\gJ^1$ is a  RAB between $\mH_{eo}$ and $\gJ^1$.

\end{theorem}

\no {\bf Proof.} i) Let $f=(S\m)\ts \gr$. From Theorem \ref{Taa}, we deduce that:

\no $\bu$ The mapping $f$ is real analytic local isomorphism.

\no $\bu$ The mapping $f$ is one-to-one.

\no $\bu$  The mapping $f$ is a surjection.

Thus we obtain that $f$ is a RAB between $\mH$ and $\gJ^o\ts \ell^2$.
From Lemma \ref{Te3} we have the identity $(\cS\m)\ts \gr=(\t \ts \gt)\circ \gF$. Then the mapping also $\t \ts \gt$ acting from $\mH$ to $\gJ^1\ts \ell^2$ is a RAB between $\mH$ to $\gJ^1\ts \ell^2$.
The proof of ii) is similar and is based on Theorem \ref{Tde}.
 \BBox


\subsection {Proof of main Theorems \ref{T1}-\ref{T5}.}

Consider inverse problems for 4-spectra mappings.

\no {\bf Proof of Theorem \ref{T1}.}
Each  $v\in\mH$ has the even-odd extension $\wt
v\in\wt\mH_{eo}$ on the interval $(0,2)$  given by \er{q}.
For the space $\wt\mH_{eo}$ we define the gap length mapping
$\wt\p_{(c)}:\wt\mH_{eo}\to \ell^2$ given by
\[
\lb{pe1}
\tes
v\to \wt\p_{(c)}=(\wt\p_{c,n})_{n\in \Z}, \qq \wt\p_{c,n}={1\/2}(\wt\l_n^-+\wt\l_n^+)-\wt\m_n,
\]
where $\wt\l_n^\pm$ are 4-periodic eigenvalues and  $\wt\m_n$ are Dirichlet
eigenvalues for the vector $\wt v\in\wt\mH_{eo}$.
By Theorem \ref{Tp2}, mapping  $v\to \wt\p_{(c)}$ is a RAB   between $\wt\mH_{eo}$ and $\ell^2$. Note that \er{GN} gives $\wt\p_{c,n}={1\/2}(\wt\n_n-\wt\m_n)$.
Then due to \er{Ev1},  the
components of the 4-spectra mapping $\gf=(\gf_n)_{n\in\Z}$ satisfy
$$
\begin{aligned}
& \gf_{2n-1}(v)=\vr_n(v) -\t_n(v)=\n_{2n-1}(\wt v)
-\m_{2n-1}(\wt v)=\wt\p_{c,2n-1}(\wt v),
\\
& \gf_{2n}(v)=\n_n(v) -\m_n(v)=\n_{2n}(\wt v) -\m_{2n}(\wt v)
=\wt\p_{c,2n}(\wt v),
\end{aligned}
$$
for all $n\in \Z$,
which yields $\gf(v)=\wt\p_c(\wt v)$. Then due to Theorem \ref{Tp2},  the mapping  $\gf: \mH \to \ell^2 $ is a RAB between the spaces
$\mH$ and $\ell^2$. Estimates \er{pe1} and the
identities $\int_0^2 \wt v^2dx=2\int_0^1 v^2dx$ and \er{esqp} yield
\er{egf}. \BBox

Now we describe properties of the mapping  $\cU_\s=\gf^{-1}\s\gf: \mH\to \mH$ for some operator $\s=(\s_j)_{n\in\Z}\in \gS$. In Section 3 we have proved some its properties. For example, from Lemma \ref{Tet} we obtain that
the   4-periodic   eigenvalues
$\{\widetilde \l_{n}^\pm,n\in\Z\}$ are invariant under $\cU_\s$
and
\begin{equation}
\label{11x}
(\widetilde \l_{n}^\pm)_{n\in\Z}=(\widetilde \l_{n}^\pm)_{n\in\Z}\circ \cU_\s.
\end{equation}

\no {\bf Proof of Theorem \ref{T2}.} i) By Theorem \ref{T1}, the mapping
$\cU_\s=\gf^{-1}\ci ( \s\gf): \mH\to \mH$ is a RAB of $\mH$ onto itself.
Due to the estimate
\er{egf} the mapping $\cU_\s=\gf^{-1}\ci ( \s\gf): \mH\to \mH$ is bounded
in any ball $\{\|v\|\le r\}$. The definition $\cU_\s=\gf^{-1} \s\gf:
\mH\to \mH$ implies that $\cU_\s=\cU_\s^{-1}$. The definition
$\cU_\s=\gf^{-1} \s\gf$ implies $\cU_\s\ci
\cU_{\s'}=\cU_{\s\s'}=\cU_{\s'}\ci \cU_{\s}$ for all $\s, \s'\in \gS.$

Due to (\ref{11x}) the   4-periodic   eigenvalues $\{\widetilde \l_{n}^\pm,n\in\Z\}$ are invariant under $\cU_\s$ and then the
Lyapunov function for the potential $\widetilde  v(x)\in\mH_{eo}$ given by \er{q} is
also invariant under $\cU_\s$. Thus the norm  $\int_0^2|\widetilde
v(x)|^2dx$ is invariant under $\cU_\s$ (see e.g., \cite{K06}) and we
obtain for $u=\cU_\s (v)$:
$$
2\int_0^1|u(x)|^2dx=\int_0^2|(\wt  u)(x)|^2dx=\int_0^2|\wt
v(x)|^2dx=2\int_0^1| v(x)|^2dx,
$$
which yields $\|\cU_\s(v)\|=\|v\|$.

\no ii)  The statement  \er{dU}  follows from \er{eig}. We show \er{peU}.
 Consider the even case $\s^e\in\gS$ when  $\s_n^e=-1$  for all odd $n\in\Z$ and $\s_n^e\in \{\pm  1\}$ for all even  $n\in\Z$. Let $ v^\bu:=\cU_{\s^e} (v)$.
 Then from Lemma \ref{Tet} for $n=2j-1, j\in\Z $ we have that
$  \tau_j(v)=\varrho_j(v^\bu),\
\varrho_j(v)=\tau_j(v^\bu)
$.
These identities and \er{FS} imply $(\l_n^\pm)_{n\in\Z}=(\l_n^\pm)_{n\in\Z}\ci \cU_{\s^e}$, since
$$
2\D(\cdot,v)=\varphi_2(1,\cdot,v)+\vartheta_1(1,\cdot,v)
=\vartheta_1(1,\cdot,q^\bu)+\varphi_2(1,\cdot,v^\bu)=2\D(\cdot,v^\bu).
$$
 Consider the odd case $\s^o\in\gS$ when  $\s_n^{o}=-1$  for all even $n\in\Z$ and $\s_n^{o}\in \{\pm  1\}$ for all  odd $n\in\Z$.
 We have the identity $\s^o =(-I)\s^e$ for some even $\s^e$. Then
 $\l_n^\pm=\l_n^\pm\ci \cU_{\s^o}$ for all $n\in\Z$, since we have the same
 for $\s^e$ and $\s=-I$ due to \er{per}.
 Finally, any $\s\in\gS$ has the form $ \s=\s^o\s^e$ for some    $\s^o, \s^e\in\gS$. Then we obtain $(\l_n^\pm)_{n\in\Z}=(\l_n^\pm)_{n\in\Z}\ci \cU_{\s}$, since $\s^o, \s^e$ keep the eigenvalues $(\l_n^\pm)_{n\in\Z}$.

\no iii) Results \er{U1}, \er{U2} follow from Lemma \ref{Tet}.
\BBox



\no {\bf Proof of Theorem \ref{T3}.} i) We show that the 2-spectra mapping
$\t\star \m: \mH\to \gJ$ is a RAB between $\mH$ and $\gJ$.
From Lemma \ref{Tet} we have the identity $\wt \m=\t\star \m$.
Then Theorem \ref{Tde} gives that the 1-spectra mapping $v\to \t\star \m$ acting from $\mH$ into $\gJ$ is a RAB   between $\mH$ and $\gJ$.
Thus from the identities \er{enm4} we obtain that
 all 2-spectra mappings $\t\star \m, \varrho\star \n$,  $\varrho\star
\m$ and $\t\star \n$ acting from $\mH$ into $\gJ$  are isomorphic, each of them  is a RAB between $\mH$ and $\gJ$ and they satisfy \er{enm4}.

ii) Due to Theorem \ref{Taa} the  mapping $\m\ts \gr $ acting from $\mH$ into $\gJ^o\ts \ell^2$ is a RAB   between $\mH$ and $\gJ^o\ts \ell^2$.
This and the identity \er{enm14} imply  \er{enm1} and the statement ii).

iii) Due to Theorem \ref{Tsa}
the mappings $\t \ts \gt$ acting from $\mH$ to $\gJ^1\ts \ell^2$ is  a RAB between $\mH$ to $\gJ^1\ts \ell^2$.
This and the identity \er{enm14x} imply  \er{enm2} and the statement iii).
 \BBox

We discuss inverse problems for 1-spectra mappings and normalizing mappings.

\begin{theorem}
\label{T???}
i) Each of the mappings $\m\ts \ga$ and $\n \ts \gb$  (defined by
\er{mntz1}, \er{mnc}) acting from $\mH$ into $\gJ^o\ts \ell^2$ is a
bijection between $\mH$ and $\gJ^o\ts \ell^2$ and they satisfy
\begin{equation}
\label{nc1}  \m\ts \ga=(\n \ts \gb)\circ \cF_o.
\end{equation}
\no ii) Each of the mappings $\t\ts \gc$ and $\vr \ts \gd$  (defined by
\er{mntz1}, \er{mnc}) acting from $\mH$ into $\gJ^1\ts \ell^2$ is a
bijection between $\mH$ and $\gJ^1\ts \ell^2$ and they  satisfy
\begin{equation}
\label{nc2}  \t\ts \gc=(\vr \ts \gd)\circ \cF_o.
\end{equation}

\end{theorem}

\no {\bf Proof.}
We show i), the proof of ii) is similar.
Consider the mappings $g:=\m\ts \ga$ and $f:=\m \ts \gr$, where
the mapping $f:\mH\to \cJ^o\ts \ell^2$ is a bijection between $\mH $ and $\cJ^o\ts\ell^2$. Due to \er{ss1} the sequence $\ga_n, n\in \Z$ satisfies $\ga_n=\gr_n-\ga_n ^\bu,$ where \er{aa1}  gives $\ga^\bu\in \ell^2$.

We show an injection. We assume that $g(v)=g(u)$ for $v, u\in\mH$.
Then we have $f(v)=f(u)$, which yields $v=u$, since $f$ is a bijection.

We show a surjection. Let $(\hat \m,\hat \ga)\in \cJ^o\ts\ell^2$.
Define the sequence $\hat \gr=\hat \ga+\hat\ga ^\bu$, where $\hat\ga^\bu
=(\m_0^o-\l) \ {\rm v.p.}{\prod_{n\in \Z,n\ne0}} {\mu_n-\l\/\mu_n^o}$
 and due to \er{aa1} it satisfies $\hat\ga^\bu\in \ell^2$.
 This gives $\hat\gr\in \ell^2$.
For $(\hat \m,\hat \gr)\in \cJ^o\ts\ell^2$ there exists $v\in\mH$
such that $f(v)=(\hat \m,\hat \gr)$. Thus we obtain $g(v)=(\hat \m,\hat \ga)$.
\BBox

\no {\bf Proof of Theorem \ref{T4}.}
By Theorem \ref{Tsa}, the mappings $(S\m)\ts \gr $ and $\t \ts \gt$ acting from $\mH$ to $\gJ^1\ts \ell^2$ are  a RAB between $\mH$ to $\gJ^1\ts \ell^2$ and satisfy $(\cS\m)\ts \gr=(\t \ts \gt)\circ \cF$.
From this and the identity \er{enm1} $
\m\ts \gr =(\n \ts \gs )\circ \cF_o$ we obtain the third identity
$(\cS\m)\ts \gr=(\t \ts \gt)\circ \cF=((\cS\n) \ts \gs )\circ \cF_o$.
Moreover, using this and the identity \er{enm2} $ \t\ts \gt =(\vr \ts \gu )\circ \cF_o$ have \er{enz1}. Then the identity \er{enz1} and the bijection
of the mappings $(S\m)\ts \gr $ from Theorem \ref{Tsa} gives
that the  mappings
$(\cS\m)\ts \gr, (\cS\n)\ts \gs, \t \ts \gt$ and $\vr \ts \gu$ acting from $\mH$ into $\gJ^1\ts \ell^2$ is a RAB between $\mH$ and $\gJ^1\ts \ell^2$.
\BBox

 \no {\bf Proof of Corollary \ref{T5}}
Due to Lemma \ref{Te4}, the mapping $\gF :\mH_{eo}\to \mH_{oe}$ is a bijection. Then the proof follows  Theorem    \ref{T4} and Lemma \ref{Te4}, iii).
 \BBox

Recall that the symplectic form has the form
$
f\wedge g=\int _0^1 (f(x),g(x))_0dx,\  f,g\in \mH
$, where $(a,b)_0=a_1b_2-a_2b_1$, for $a,b\in\C^2$.
We show  the following canonical relations.

\begin{theorem}  \lb{Tcv}
 For any $n, j\in \Z$ and $v\in\mH$ the following identities  hold true:
\[
\lb{cv1}
\begin{aligned}
\n_n'(v)\we \n_j'(v)=0,\qqq
\gs'_{n}(v)\we \n_j'(v)=\d_{nj},
\qqq
\gs'_{n}(v)\we \gs'_{j}(v)=0,
\\
\t_n'(v)\we \t_j'(v)=0,\qqq
\gt'_{n}(v)\we \t_j'(v)=\d_{nj},
\qqq
\gt'_{n}(v)\we \gt'_{j}(v)=0,
\\
\vr_n'(v)\we \vr_j'(v)=0,\qqq
\gu'_{n}(v)\we \vr_j'(v)=\d_{nj},
\qqq
\gu'_{n}(v)\we \gu'_{j}(v)=0,
\end{aligned}
\]
where $\n_n'={\pa\n_n\/\pa v},\gs'_{n}={\pa \gs_n \/\pa v}$ and each of  sequences $\{ \n_n', \gs'_{n}\}_{n\in \Z},  \{ \t_n', \gt'_{n}\}_{n\in \Z},.. $,
is a basis for $\mH$.
\end{theorem}

\no {\bf Proof.} Let we have functions $a:\mH\to \R$ and $b:\mH\to \R$
and $a(v), b(v) v\in\mH$.
Consider the linear mapping $F:\mH\to\mH$, given by $Fv=F(x)v(x)$, where
$F(x)$ is some $2\ts 2$ matrix. Define functions $A(u)=a(Fu)$ and $B(u)=b(Fu)$  for $u\in\mH$. We have $a_v'=(a_{v_1}',a_{v_2}')$, where $a_{v_j}'={\pa\/\pa v_j}a, j=1,2$. Then we obtain $A'(u)\we B'(u)=\int _0^1 (A_u',B_u')_0dx$, where
\[
\lb{FSz}
\begin{aligned}
 (A_u',B_u')_0= (A_{u_1}',A_{u_2}')  J (B_{u_1}',B_{u_2}')^\top=
(a_{v_1}',a_{v_2}') F J F^\top (b_{v_1}',b_{v_2}')^\top,
\end{aligned}
\]
We show the first line in \er{cv1}. Due to \er{enm4} we have
$\n=\m\circ \cF_o$ and $\gs= \gr\circ \cF_o$.
Then using \er{mr}, \er{FSz} and  $\cF_o J \cF_o^\top=J$  we obtain
$$
\begin{aligned}
 \n_n'(u)\we \n_j'(u)=\m_n'(v)\we \m_j'(v)=0,
 \qqq
 \gs_n'(u)\we \gs_j'(u)=\gr_n'(v)\we \gr_j'(v)=0,
 \\
 \gs_n'(u)\we \n_j'(u)=\gr_n'(v)\we \m_j'(v)=\d_{n,j}.
 \end{aligned}
$$
We show the second line in \er{cv1}. Due to \er{cA} we have
$\t=\cS \m\circ \gF$ and $\gt= \gr\circ \gF$.
Then using \er{mr}, \er{FS} and  $\gF J \gF^\top=J$  we obtain
$$
\begin{aligned}
 \t_n'(u)\we \t_j'(u)=\m_n'(v)\we \m_j'(v)=0,\qq
 \gt_n'(u)\we \gt_j'(u)=\gr_n'(v)\we \gr_j'(v)=0,
 \\
 \gt_n'(u)\we \m_j'(u)=\gr_n'(v)\we \m_j'(v)=\d_{n,j}.
 \end{aligned}
$$
The identity \er{enm14x} gives $\t\ts \gt =(\vr \ts \gu )\circ \cF_o$. Then
using similar arguments we obtain the last line in \er{cv1}.
The proof  about the basis repeats the case of $\m_n', \gr_n', n\in \Z$
from \cite{K01}.
\BBox

We discuss asymptotics of spectral data.

\begin{theorem}
\label{Tas}
For beach $d\in(1,2)$.
The mappings $\n\ts \gs$, $\t\ts \gt$, $\vr\ts \gu$ have
following asymptotics
\[
\lb{ass1}
\begin{aligned}
\ma \n_n(v)-\pi n \\ \gs_n(v) \am =J_1(\F v)_n+\ell^d(n),
\end{aligned}
\]
\[
\lb{ass2}
\begin{aligned}
 \ma \t_n(v)-\t_n^o \\ \gt_n(v) \am=-J_1(\F\gF^* v)_n+\ell^d(n),
\end{aligned}
\]
\[
\lb{ass3}
 \ma \vr_n(v)-\vr_n^o \\ \gu_n(v) \am=J_1(\F\gF^* v)_n+\ell^d(n),
\]
as $n\to\pm \iy$, uniformly on $\cB_\C(u,\ve_{u})$, for any $u\in\mH$
 and $\ve_{u}=4^{-4}e^{-3\|u\|}$.

\end{theorem}

\no {\bf Proof.}
Theorem \ref{T3} gives the identity  $\n\ts\gs=(\m\ts\gr)\ci \cF_o $. Then \er{ra1} yields that
$$
\ma \n_n (v)-\pi n\\ \gs_n (v)\am=
\ma \m_n (w)-\pi n\\ \gr_n (w)\am= -J_1(\F w)_n+\ell^d(n),\qqq v=\cF_o w,
$$
uniformly on $w\in \cB_\C(q,\ve_{q}), q\in\mH$. This implies \er{ass1}
since $\|v-u\|=\|w-q\|$, where $u=\cF_oq$.

Theorem \ref{T3} gives that $(\t\ts\gt)\circ\gF=(\cS\m)\ts\gr$. Then \er{ra1} yields that
$$
\begin{aligned}
\ma \t_n (v)\\ \gt_n (v)\am=\ma \cS\m_n (w)\\ \gr_n (w)\am
=\ma \t_n^o \\ 0 \am -J_1(\F w)_n+\ell^2(n), \qqq v=\gF^* w,
\end{aligned}
$$
uniformly on $w\in \cB_\C(q,\ve_{q}), q\in\mH$. This implies \er{ass2}
since $\|v-u\|=\|w-q\|$, where $u=\gF^*q$.

 Theorem \ref{T3} gives the identity $(\vr\ts\gu)\circ\cF_o=\t\ts\gt$.
 Then asymptotics of  $\t\ts\gt$ and similar arguments imply \er{ass3}.
  \BBox


\section {Isomorphic inverse problems on the circle \lb{Sec5}}
\setcounter{equation}{0}


\subsection{Periodic potentials} We prove the first results about periodic inverse problems. We apply results from Theorem \ref{T1} to the periodic inverse problems.


\begin{proposition}
\lb{Tpnb}
 Let  eigenvalues $(\l_{2n}^\pm(v))_{n\in\Z}$ and one of the following be given
 for some $v\in\mH$:

\no  i) $\t_n(v)$ and  $\sign (\m_n(v)-\n_n(v))$ for all
$n\in \Z$.

\no  ii) $\vr_n(v)$ and  $\sign (\m_n(v)-\n_n(v))$ for all
$n\in \Z$.

\no  iii) $\m_n(v)$ and  $\sign \ln |\vp_2(1,\m_n(v),v)|$ for all
$n\in \Z$.

\no  iv) $\n_n(v)$ and  $\sign \ln |\vt_1(1,\n_n(v),v)|$ for all
$n\in \Z$.

Then  the potential $v$ is uniquely determined.

\end{proposition}

 \no {\bf Proof}.  i) Let $v\in\mH$.
 It is known that the function $\D(\l,v)-1$ is recovered by its zeros, i.e., the periodic spectrum $\l_{2n}^\pm(v), n\in \Z$.   Due to \er{FS} we have
$\vp_2(1,\l)=  {\rm v.p.}\prod_{n\in \Z}{\t_n-\l\/\t_n^o}$.
Thus using the definition of the Lyapunov function $\D(\l,v)=
{1\/2}(\vp_2(1,\l,v)+\vt_1(1,\l,v))$, we can recover the function $\vt_1(1,\l,v))$ and its zeros $\vr_n(v), n\in \Z$.

Let $\wt v=\cE v\in \wt\mH_{eo}$ be an even extension of $v$
given by \er{q} and let
$\wt \l_n^\pm$  be corresponding  periodic eigenvalues.
From Femma \ref{Tet} we have
$$
\{\wt \l_{2n-1}^-, \wt \l_{2n-1}^+\}=\{\vr_{n},\t_{n}\}, \qq
\{\wt \l_{2n}^-,\wt \l_{2n}^+\} =\{\m_{n},\n_{n}\},\qq \forall \ n\in \Z.
$$
Thus using the anti-periodic eigenvalues   $\{\wt \l_{2n-1}^-, \wt \l_{2n-1}^+\}=\{\vr_{n},\t_{n}\} $  we determine the periodic eigenvalues $\{\wt \l_{2n}^-,\wt \l_{2n}^+\} =\{\m_{n},\n_{n}\}$ for $\wt v$.
This jointly with the sequence $\sign (\m_n(v)-\n_n(v)),  n\in \Z$ gives $\m_n(v), \n_n(v)$ for all $ n\in \Z$. Moreover, due to Theorem \ref{T3}
the potential $v$ is uniquely determined. The proof of ii)-iv) is simlar
and iv) is well known, see e.g., Theorem \ref{TKE}.
 $\BBox$


We consider inverse problems on the circle. Firstly, we define the
gap mapping $v\to \p=(\p_n)_{n\in \Z}$ acting from $\mH$ into $ \ell^2\os
\ell^2$ from \cite{K99}. The components $\p_n\in \R^2$  are
constructed via the periodic plus Dirichlet eigenvalues plus signs
by
\begin{equation}
\label{gL1}
\begin{aligned}
&  \textstyle   \p_n=(\p_{c,n},\p_{s,n})\in\R^2, \qqq
|\p_n|^2=\p_{c,n}^2+\p_{s,n}^2= {1\/4}(\l_n^+-\l_n^-)^2,
\\
&     \textstyle \p_{c,n}={1\/2}(\l_n^++\l_n^-)-\mu_n,\qq
\p_{s,n}=\big||\p_n|^2- \p_{c,n}^2  \big|^{1\/2}\sign \gr_n ,\qq
\gr_n =\log|\vp_2(1,\mu_n)|.
\end{aligned}
\end{equation}
The mapping $\p$ is a RAB between $\mH $ and $\ell^2\os \ell^2$, see
Theorem \ref{TKE} below.

 We define another gap mapping   $\gp: \mH \to \ell^2\os
\ell^2$  by $v\to \gp=(\gp_n)_{n\in \Z}$. The components $\gp_n\in \R^2$
are constructed via the 2-periodic $\l_n^\pm$ plus Neumann eigenvalues $\n_n$ plus $\sign \gs_n$ by
\begin{equation}
\label{gL2}
\begin{aligned}
&\textstyle  \gp_n=(\gp_{c,n},\gp_{s,n})\in\R^2,\qq
|\gp_n|^2=\gp_{c,n}^2+\gp_{s,n}^2= {1\/4}(\l_n^+-\l_n^-)^2,
\\
& \textstyle \gp_{c,n}={1\/2}(\l_n^++\l_n^-)-\nu_n,\qqq
\gp_{s,n}=\Big||\gp_n|^2- (\gp_{c,n})^2  \Big|^{1\/2}\sign
\gs_n ,\qq \gs_n =\ln |\vartheta_1(1,\nu_n)|.
\end{aligned}
\end{equation}

Secondly we consider inverse problems in terms of local maxima and
minima of the Lyapunov function, given by
$\D(\l)={1\/2}(\vp'(1,\l)+\vt(1,\l))$. The Lyapunov function on the
real line  has local maxima and minima at points $\l_n\in
[\l_n^-,\l_n^+]$ for all $n\in \Z$, where $(-1)^n\D(\l_n^\pm)=1$ and
$(-1)^n\D(\l_n)\ge 1$. Define the corresponding mapping $h: \mH\to
\ell^2\os\ell^2$ as $h: v\to h=(h_n)_{n\in \Z}$  from
\cite{K01}. The components $h_n=(h_{c,n}, h_{s,n})\in \R^2$ are constructed via maxima
and minima of the Lyapunov function plus Dirichlet eigenvalues plus
signs by
\begin{equation}
\label{h75}
\begin{aligned}
 h_{c,n}=\Big||h_n|^2-h_{s,n}^2\Big|^{1\/2}{\rm sign} (\l_n-\mu_n),
 \qqq h_{s,n}=\gr_n=-\log |\vp_2(1,\m_n)|.
 \end{aligned}
\end{equation}
The value $|h_n|^2=h_{c,n}^2+h_{s,n}^2\ge 0$ is uniquely defined by the equation
$\ch |h_n|=|\D(\l_n)|\ge 1$. Recall that $(-1)^{n}\D(\mu_n)=\ch
h_{s,n}$ for all $n\in\Z$ and $|h_n|\ge |h_{s,n}|$, since
$(-1)^n\D$ has the maximum at $\l_n$ on the segment
$[\l_n^-,\l_n^+]$. The mapping $h$ is a RAB between $\mH$ and
$\ell^2\os \ell^2$.

 We introduce similar mapping $\gh: \cH\to \ell^2\os\ell^2$ as $\gh: v\to
\gh(v)=(\gh_n(v))_{n\in \Z}$. The components $\gh_n=(\gh_{c,n}, \gh_{s,n})\in
\R^2$  are constructed via maxima and minima of the Lyapunov
function plus Neumann eigenvalues plus signs by
\begin{equation}
 \label{gh75}
\begin{aligned}
\gh_{c,n}=\Big||\gh_n|^2-\gh_{s,n}^2\Big|^{1\/2}{\rm sign}
(\l_n-\nu_n),
 \qqq \gh_{s,n}=- \log |\vartheta_1(1,\nu_n)|.
 \end{aligned}
\end{equation}
Recall that $(-1)^{n}\D(\nu_n)=\ch \gh_{s,n}$ for all $n\in\Z$ and
$|\gh_n|\ge |\gh_{s,n}|$, since $(-1)^n\D$ has the local maximum at
$\l_n$ on the segment $[\l_n^-,\l_n^+]$. Recall results from \cite{K01}, \cite{K05}.

\begin{theorem}
\label{TKE} \no i) The mapping $h: \mH  \to
\ell^2\os \ell^2$ given by (\ref{h75}) is a RAB between $\mH $ and
$\ell^2\os \ell^2$. Furthermore, the following estimates hold true:
\begin{equation}
\label{esqh}
\tes {1\/2}\|v\| \le \|h\| \le
3\|v\|(1+\|v\|)^{1\/2},
\end{equation}
where $\|v\|^2=\int_0^1v^2(x)dx$ and  $\|h\|^2=\sum_{n\in \Z}|
h_{n}|^2$.

\no ii) The mapping $\p: \mH  \to \ell^2\os \ell^2$ given by (\ref{gL1})
is a RAB between $\mH $ and $\ell^2\os \ell^2$. Furthermore, the
following estimates hold true:
\begin{equation}
\label{esqp}
\tes
{1\/\sqrt2}\|\p\|\le  \|v\| \le 2\|\p\|(1+\|\p\|),
\end{equation}
where  $\|\p\|^2=\sum_{n\in \Z}
(\p_{c,n}^2+\p_{s,n}^2)={1\/4}\sum_{n\in \Z}|\l_n^+-\l_n^-|^2$.

\no iii) Let $v\in \mH$ and $
Q_2=\int_0^1(|v'|^2+|v|^4)dx$ and $\|\p\|_1^2=\|\p\|^2+\sum_{n\in \Z}(2\pi n)^2|\p_{n}|^2
$. Then
\[
\label{esqp1}
\tes
{1\/24}\|\p\|_1^2\le Q_2\le 8\Big((\pi+\|v\|^2) \|\p \|_1^2+\|v\|^2\Big).
\]
\end{theorem}

We describe isomorphic mappings on the circle.

\begin{corollary}
\label{Tp1}

\no i) The mappings $\p:\mH\to \ell^2\os \ell^2$ and
$\gp:\mH\to \ell^2\os \ell^2$   are isomorphic, each of them  is a RAB between $\mH$ and $\ell^2\os \ell^2$ and they satisfy
\begin{equation}
\label{pf1}
\begin{aligned}
\p=\gp\circ \cF_o.
\end{aligned}
\end{equation}
\no ii)  The mappings $h:\mH\to \ell^2\os \ell^2$ and
$\gh:\mH\to \ell^2\os \ell^2$   are isomorphic, each of them  is a RAB between $\mH$ and $\ell^2\os \ell^2$ and they satisfy
\begin{equation}
\label{pf2}
\begin{aligned}
h=\gh\circ \cF_o.
\end{aligned}
\end{equation}

\end{corollary}
\no {\bf Proof.} The proof follows from Theorem \ref{TKE} and Lemma \ref{Te1}.
\BBox

 We formulate the key result of the direct method, proved in \cite{KK97},
incorporating a necessary modification from \cite{K01}.


\begin{theorem} \label{Tfa2}
Let $H, H_1$ be real separable Hilbert spaces equipped with norms
$\|\cdot \|, \|\cdot \|_1$ respectively. Suppose that a map $f: H \to H_1$
satisfies the following conditions:

\no i) $f$ is real analytic,

\no ii) the derivtive $f'$ has an inverse for all $v\in H$,

\no iii) there is a nondecreasing function
$\x : [0, \iy ) \to [0, \iy ), \x (0)=0,$ such
that $\|v\|\leq \x (\|f(v)\|_1)$  for all $v\in H$,

\no iv) there exists a basis $\{e_n\}_{n\in \Z}$ of $H_1$ such that
each map  $(f(\cdot ), e_n)_1: H \to \R, n\in \Z$,  is compact,

\no v) for each $C>0$ the set $\{v\in H: \sum_{n\in \Z} n^2(f(v), e_n)_1^2<C \}$
is compact.

\no Then $f$ is a real analytic isomorphism between $H$ and $H_1$.
\end{theorem}

\begin{theorem}
\label{Tp2}

A gap lenght mapping $\p_{(c)}:\mH_{eo}\to \ell^2$
given by
\[
\lb{pe1w}
\tes
v\to \p_{(c)}=(\p_{c,n})_{n\in \Z}, \qq \p_{c,n}={1\/2}(\l_n^-+\l_n^+)-\m_n,
\]
is a RAB   between $\mH_{eo}$ and $\ell^2$.

\end{theorem}
\no {\bf Proof.} In order to prove theorem we use Theorem \ref{Tfa2} and
check all its condiotions.

i) In Theorem \ref{TKE} we proved that the mapping $v\to \p_{(c)}$ is real
analytic.

ii) the derivtive $f'$ has an inverse for all $v\in\mH_{eo}$,

iii) We have the needed estimates $\|v\| \le 2\|\p_{(c)}\|(1+\|\p_{(c)}\|)$ follows from \er{esqp}.

iv)  Using Lemma \er{GN} we have the identity
$\p_{c,n}={1\/2}(\l_n^-+\l_n^+)-\m_n={1\/2}(\n_n+\m_n)$ for all $n\in\Z$.
Lemma \ref{T23} gives that each
map  $\m_n: \mH_{eo} \to \R, n\in \Z$  is compact.

v) Using  \er{esqp1} we deguce that the set $\{v\in \mH: \|\p(v)\|_1\le C   \}$ is compact for each $C>0$.

Then by Theorem \ref{Tfa2}, $f$ is a real analytic isomorphism between $\mH_{eo}$ and $\ell^2$.
\BBox

\no {\bf{Proof of Corollary  \ref{Ter1}}.} i) Recall that
$\cU_\s=\gf^{-1}\circ(\s\gf)$ is defined by \er{dj}.
From Lemma \ref{Te1}, and \er{Rep}, \er{Rep1} we deduce that
$$
\t_n\circ\cU_\s =\ca \vr_n, \ \ n\in\Y_1, \\
\t_n, \ \ n\notin\Y_1, \ac \qqq
\m_n\circ\cU_\s =\ca \n_n, \ \ n\in\Y_2, \\
\m_n, \ \ n\notin\Y_2 \ac.
$$
These identities give \er{rtm}.
Theorem \ref{T3} shows that the mapping $\t\star\mu$ is a RAB, then the mapping
$\z\star \phi:\mH \to \gJ$ is a RAB between $\mH$ and $\gJ$.
 \BBox




\begin{thebibliography} {9999}\setlength{\itemsep}{-\parskip}
\footnotesize \footnotesize


\bibitem {AT04} M. Ablowitz; B. Prinari, A.
Trubatch, Discrete and continuous non-linear Schr\"odinger systems,
2004.

\bibitem {AHM05}
Albeverio, S.; Hryniv, R.; Mykytyuk, Ya. Inverse spectral problems for Dirac operators with summable potentials. Russ. J. Math. Phys. 12 (2005), no. 4, 406--423.

\bibitem {A14} Amour, L.
The coordinate system $\m\ts k$ on $L^2([0,1])\ts L^2([0,1])$  for the AKNS operator, Math. Phys. Anal. Geom. 17(2014), 83--93.

\bibitem {A93}
Amour, L. Inverse spectral theory for the AKNS system with separated boundary conditions. Inv. Probl. 5(1993), 507--523.

\bibitem {BG93} D. Battig, B. Grebert, J. Guillot, T. Kappeler:
Foliation of phase space for the cubic non-linear  Schr\"odinger
equation, Compositio Math., 1993, 85, 163-199.


\bibitem {CKK04} Chelkak, D.; Kargaev, P.; Korotyaev, E. Inverse problem
for harmonic oscillator perturbed by potential,
characterization. Comm. Math. Phys. 249 (2004), no. 1, 133-196.


\bibitem {CK07}    Chelkak,  D.; Korotyaev, E. The inverse problem for
perturbed harmonic oscillator on the half-line with Dirichlet
boundary conditions, Ann. Henri Poincare 8(2007), no.6, 1115--1150.

\bibitem {CK06}
  Chelkak,  D.; Korotyaev, E. Parametrization of the isospectral set
for the vector-valued Sturm-Liouville problem, J. Funct. Anal.
241(2006), 359--373.


\bibitem {CK09} Chelkak, D.; Korotyaev, E.  Weyl-Titchmarsh
functions of vector-valued Sturm-Liouville operators on the unit
interval. Journal of Functional Analysis, 257 (2009), 1546-1588.

\bibitem {CG02} S. Clark and F. Gesztesy, Weyl-Titchmarsh "M-Function Asymptotics, Local Uniqueness Results, Trace Formulas, and Borg-Type Theorems for Dirac Operators," Trans. Amer. Math. Soc. 354 (9)(2002), 3475--353.

\bibitem {CM17} Clay, M.; Margalit, D. Office Hours with a Geometric
Group Theorist. Princeton Univ. Press, 2017.

\bibitem{CM93} Coleman C.; McLaughlin, J. Solution of
the inverse problem for an impedance with integrable derivative, I
Commun. Pure Appl. Math. 46 (1993), 145--184.

\bibitem{CM93a} Coleman C.; McLaughlin, J. Solution of
the inverse problem for an impedance with integrable derivative II,
Commun. Pure Appl. Math. 46(1993), 185--212.


\bibitem {DK00}  Daskalov V. B.;  Khristov, E. Kh. Explicit Formulae for the Inverse Problem for t,he Regular Dirac Operator, Inverse Problems 16 (1)(2000), 247--258.

\bibitem {DG01}    del Rio, R., Grebert, B. Inverse spectral results for the AKNS systems with partial information on the
potentials.Math. Phys. Anal. Geom. 4(3)(2001), 229--244.

\bibitem{FM76} Flaschka H.; McLaughlin D. Canonically
conjugate variables for the Korteveg- de Vries equation and the Toda
lattice with periodic boundary conditions. Prog. of Theor. Phys.
55(1976),  438--456.


\bibitem {GT87} Garnett, J.; Trubowitz, E. Gaps and bands of one dimensional
periodic Schr\"odinger operators II. Comment. Math. Helv. 62(1987),
18--37.

\bibitem {GD75} Gasymov, M.G.;  Dzhabiev T. T.  Determination of a Dirac system of differential equations from two spectra,  Trudy Letnei Shkoly po Spectral'noi Teorii Operatorov i Teorii Predstavleniya Group, ELM Baku, 1975, 46--71.

\bibitem {GL66}
Gasymov, M. G.;  Levitan, B. M. The inverse problem for the Dirac system, Dokl. Akad. Nauk SSSR, 167:5 (1966), 967--970.

\bibitem {GL51} Gel'fand, I. M.; Levitan, B. M. On the determination of a differential equation from its spectral function. Izvestiya Akad. Nauk SSSR. Ser. Mat. 15, (1951). 309--360.

\bibitem {GL53} Gel'fand, I. M.; Levitan, B. M. On a simple identity for the characteristic values of a differential operator of the second order. Doklady Akad. Nauk SSSR (N.S.) 88, (1953). 593--596.

\bibitem {Ge02}  Gesztesy, F.; Kiselev, A.;  Makarov, K. A. Uniqueness results for matrix-valued Schr\"odinger, Jacobi, and Dirac-type Operators, Math. Nachr. 239/240(2002), 103--145.

\bibitem {GG93}    Grebert, B.; Guillot, J. Gaps of one-dimensional periodic AKNS systems. Forum Math. 5 (1993), no. 5, 459--504.




\bibitem {GK14} B. Grebert; T. Kappeler, The defocusing NLS equation and its normal form, 2014, European Math. Soc.


    \bibitem {H09} Hilbert, D. Zur theory der korfomen abbildung,
Nach. Kgl. Ges. G\"ottingen,, Math.-phys. Ki. 1909, 314--323.


\bibitem {H01}     M. Horvath, On the inverse spectral theory of Schr\"odinger and Dirac operators, Trans. Amer. Math. Soc. 353(10)(2001), 4155--4171.

\bibitem {IT83}
Isaacson, E. L.; Trubowitz, E. The inverse Sturm-Liouville problem.
 I. Comm. Pure Appl. Math. 36 (1983), no. 6,   767--783.

\bibitem {IMT84} Isaacson, E. L.; McKean, H. P.; Trubowitz, E.
The inverse Sturm-Liouville problem. II. Comm. Pure Appl. Math. 37
(1984), no. 1, 1--11.



\bibitem {KK95} Kargaev, P.; Korotyaev, E. Effective masses and conformal
mappings. Comm. Math. Phys. 169 (1995), no. 3, 597--625.


\bibitem {KK97} Kargaev, P. ;  Korotyaev E.
The inverse problem for the Hill operator, direct approach. Invent.
Math.  129(1997), 567--593.


\bibitem {K25}
Korotyaev, E. Isomorphic inverse problems, Russian
Journal of Math. Phys., 32(2025),  No. 2, 314--340.

\bibitem {K19} Korotyaev, E. Inverse Sturm-Liouville problems for
non-Borg conditions,  Journal of Inverse and Ill-Posed Problems,
27(2019), no 3. 445--452.

\bibitem {K06}  Korotyaev, E. Estimates for the Hill operator. II. J.
Differential Equations 223 (2006), no. 2, 229--260.

\bibitem{K05} E. Korotyaev.  Inverse problem and estimates for periodic
Zakharov-Shabat systems, J. Reine Angew. Math. 583(2005), 87--115.

\bibitem {K01}
E. Korotyaev,  Marchenko-Ostrovki mapping for periodic
Zakharov-Shabat systems,   J. Differential Equations, 175(2001), no.
2, 244--274.

\bibitem {K00}   Korotyaev, E. Estimates for the Hill operator. I. J.
Differential Equations 162 (2000), no. 1, 1--26.


\bibitem {K99}  Korotyaev, E. The inverse problem and trace
formula for the Hill operator, II.  Math. Z. 231(1999), 345--368.

\bibitem {K98} Korotyaev, E. Estimates of periodic potentials in terms
of gap lengths. Comm. Math. Phys. 197 (1998), no. 3, 521--526.


\bibitem {K97} Korotyaev, E. The inverse problem for the Hill operator.
I. Int. Math. Res. Notices 1997, no. 3, 113--125.

\bibitem {K97m}  Korotyaev, E. The estimates of periodic potentials in
terms of effective masses. Comm. Math. Phys. 183 (1997), no. 2,
383--400.

\bibitem {K96}
E. Korotyaev,  Metric properties of conformal mappings on the
complex plane with parallel slits, Internat. Math. Res. Notices,
10(1996), 493--503.


\bibitem {KC09}  Korotyaev, E.; Chelkak,  D.
 The inverse Sturm-Liouville problem with mixed boundary conditions,
St. Petersburg Math.  Journal. 21(2009), no 5, 114--137.

\bibitem {KKL} Korotyaev, E.;  Leonova, E.
Isomorphic inverse problems for Jacobi operators, preprint, 2025.


\bibitem {KM22} Korotyaev, E.; Mokeev, D. Dubrovin equation for periodic
Dirac operator on the half-line, Appl. Anal. 101 (2022), No 1,
337--365.

\bibitem {Kr51}
Krein, M. G. Solution of the inverse Sturm-Liouville problem.
 Doklady Akad. Nauk SSSR (N.S.) 76, (1951). 21--24.

\bibitem {Kr54}
Krein, M. G. On a method of effective solution of an inverse
boundary problem.  Doklady Akad. Nauk SSSR (N.S.) 94,
(1954). 987--990.

\bibitem  {Le49} Levinson, N. The inverse Sturm-Liouville problem.
Mat. Tidsskr. B. 1949, (1949). 25--30.


\bibitem{L87} Levitan, B. Inverse Sturm-Liouville problems.
Utrecht: VNU Science Press, 1987.

\bibitem {LG64}  Levitan, B. M.; Gasymov, M. G. Determination of a differential equation
by two spectra. (Russian) Uspehi Mat. Nauk 19 1964 no. 2 (116), 3--63.




\bibitem{LS91} Levitan,  B. M.;  Sargsjan, I. S. Sturm-Liouville and Dirac Operators (Nauka, Moscow, 1988; Kluwer Acad., Dordrecht, 1991.

\bibitem {MT81} McKean, H. P.; Trubowitz, E. The spectral class of
the quantum-mechanical harmonic oscillator.
Comm. Math. Phys. 82 (1981/82), no. 4, 471--495.



\bibitem {M50}  Marchenko, V. A. Concerning the theory of a differential operator of the second order. Doklady Akad. Nauk SSSR. (N.S.) 72,
(1950). 457-460.


\bibitem {MO75} Marchenko, V.; Ostrovski, I. A
characterization of the spectrum
 of the Hill operator. Mat. Sb. 97(139), (1975), 540--606.

 \bibitem{M86}
 Marchenko, V. Sturm-Liouville operator and applications.
Basel, Birkh\"auser 1986.


\bibitem{Mi78} Misura T. Properties of the spectra of periodic
and anti-periodic
boundary value problems generated by Dirac operators. I, Theor. Funktsii
Funktsional. Anal. i Prilozhen, (Russian), 30 (1978), 90--10.

\bibitem{Mi79} Misura T. Properties of the spectra of periodic
and anti-periodic
boundary value problems generated by Dirac operators, II, Theor. Funktsii
Funktsional. Anal. i Prilozhen, (Russian), 31 (1979), 102--109.

\bibitem{Mi80} Misura T. Finite-zone Dirac operators.  Theor. Funktsii
Funktsional. Anal. i Prilozhen, (Russian), 33 (1980), 107--111.

\bibitem{MT02} Mochizuki K.; Trooshin, I. Inverse problem for interior spectral data of the Dirac operator
on a finite interval, Publ. Res. Inst. Math. Sci. 38 (2002), no. 2, 387--395.


\bibitem{PT87} P\"oschel, P., Trubowitz E. Inverse Spectral Theory.
Boston, Academic Press, 1987.


\bibitem{S70} I. S. Sargsjan, A uniqueness theorem for the solution of the inverse problem for a one-dimensional Dirac system, in Some Boundary Value Problems of Ordinary Differential Equations (Univ. Druzhby Narodov, Moscow, 1970), pp. 3--13.

    \bibitem{ZS72} V. E. Zakharov; A. B. Shabat, Exact theory
of two-dimensional self-focusing and one-dimensional self-modulation
of waves in nonlinear media, Sov. Phys. JETP 34 (1972), 62--69.

\end{thebibliography}
\end{document}